\documentclass[11pt]{amsart}

\author[P.~Miska]{Piotr Miska}
\address{Institute of Mathematics\\
        Faculty of Mathematics and Computer Science\\
        Jagiellonian University in Krak\'ow\\
        ul. {\L}ojasiewicza 6, 30-348 Krak\'ow\\
        Poland\\
        and
        Department of Mathematics\\
	J. Selye University\\
	P. O. Box 54\\945 01 Kom\'arno\\
	Slovakia}
\email{piotr.miska@uj.edu.pl, miskap@ujs.sk}

\author[N.~Murru]{Nadir Murru}
\address{Department of Mathematics, Università degli Studi di Trento}
\email{nadir.murru@unitn.it}

\author[G. Romeo]{Giuliano Romeo}
\address{Department of Mathematical Sciences "Giuseppe Luigi Lagrange", Politecnico di Torino}
\email{giuliano.romeo@polito.it}

\subjclass[2010]{}

\title{On the arithmetic of multidimensional continued fractions}

\usepackage{verbatim}
\usepackage{amsmath}
\usepackage{amssymb}
\usepackage{amsthm}
\usepackage{geometry}
\usepackage{hyperref}
\usepackage{enumerate}
\hypersetup{colorlinks=true}
\usepackage{todonotes}
\usepackage[utf8]{inputenc}
\usepackage[ruled,linesnumbered]{algorithm2e}

\newcommand{\N}{\mathbb{N}}
\newcommand{\Z}{\mathbb{Z}}
\newcommand{\Q}{\mathbb{Q}}

\newcommand{\R}{\mathbb{R}}

\newcommand{\sgn}{\operatorname{sgn}}

\newtheorem{thm}{Theorem}[section]

\theoremstyle{remark}

\newtheorem{Example}[thm]{Example}

\begin{document}

\begin{abstract}
The problem of developing an arithmetic for continued fractions (in order to perform, e.g., sums and products) does not have a straightforward solution and has been addressed by several authors. In 1972, Gosper provided an algorithm to solve this problem. In this paper, we extend this approach in order to develop an arithmetic for multidimensional continued fractions (MCFs). First, we define the M\"obius transform of an MCF and we provide an algorithm to obtain its expansion. Similarly, we deal with the bilinear transformation of MCFs, which covers as a special case the problem of summing or multiplying two MCFs. Finally, some experiments are performed in order to study the behavior of the algorithms.
\end{abstract}


\maketitle

\section{Introduction}\label{section1}

Given a real number $x_0$, its continued fraction expansion is
\begin{equation} \label{cf}
x_0=a_{0}+\cfrac{1}{a_{1}+\cfrac{1}{a_{2}+\cfrac{1}{a_{3}+\cfrac{1}{\ddots}}}}=[a_{0};a_{1},a_{2},a_3,\ldots],
\end{equation}
where the \emph{partial quotients} $a_k$ are obtained by 
$$
a_{k}=\lfloor x_{k}\rfloor,\quad  x_{k+1}=\frac{1}{x_{k}-a_{k}},
$$
for $k = 0, 1, \ldots$, if the \emph{complete quotient} $x_k$ is not an integer. The expansion is finite if there exists a $k$ such that $x_{k} \in \mathbb N$ and in this case $x_0 = [a_0;a_1, \ldots, a_k]$. Continued fractions are very important and useful tools, with many applications like, e.g., in Diophantine approximation and it is well--known that they provide the best approximations of real numbers. The continued fraction expansion is another way to represent real numbers. The continued fraction representation has many advantages compared to the decimal expansion or the expansions in other bases. In fact, for instance, the continued fraction expansion is finite for rational numbers and eventually periodic for quadratic irrationals. However, we usually represent real numbers through their decimal expansion because it leads to a simple arithmetic and it allows to easily perform computations among them.
Hence, it naturally arises the problem of performing sums and products of real numbers by means of their continued fraction expansion, i.e., to develop an arithmetic of continued fractions. In general, given two continued fractions $[a_0; a_1, \ldots]$ and $[b_0; b_1, \ldots]$, finding the continued fraction expansion of
\[ [a_0; a_1, \ldots] + [b_0; b_1, \ldots], \quad  [a_0; a_1, \ldots] \cdot [b_0; b_1, \ldots]\]
is not straightforward.
This problem has been addressed by several authors. The continued fraction expansion of $2x_0$ has been investigated by Hurwitz \cite{Hur} in 1891. After that, one of the first works about the arithmetic of continued fractions is due to Hall \cite{Hall}, where the author proved that every real number can be expressed as a sum of two continued fractions with partial quotients less than or equal to $4$ (except for the first partial quotient) and similarly for the product. Then, the author provided a first algorithm for obtaining the continued fraction expansion of the linear fractional transformation of a given continued fraction (this allows, e.g., to compute the product or the sum of a continued fractions by a rational). Finally, Hall applied this algorithm for obtaining further results about sum and products of continued fractions. However, these algorithms are not straightforward and hard to use in practice. Raney \cite{Raney} proposed an alternative approach for studying the linear fractional transformation of a continued fraction based on the $(L, R)$-sequence expansion of a continued fraction. In \cite{Gosper}, Gosper proposed a straightforward algorithm for dealing also with bilinear fractional transformation of two continued fractions (sum and product of two continued fractions is a special case). Further related works are \cite{Col, Cusick, LiarStam}.

The notion of continued fraction has been generalized to higher dimensions by Jacobi \cite{Jacobi} and Perron \cite{Perron}, aiming at solving the Hermite problem \cite{HER}, i.e., to provide a generalization of the continued fraction algorithm such that it becomes eventually periodic if and only if the input is an algebraic irrational of any degree (as well as it happens for classical continued fractions and quadratic irrationals). 
Namely, given an $m$-tuple $(x^{(1)}_0,...,x^{(m)}_0)$ of real numbers, its multidimensional continued fraction (MCF) expansion is the $m$-tuple $\textbf{a}=\left[\textbf{a}^{(i)}\right]_{i=1}^m=[(a_n^{(1)})_{n\in\N},...,(a_n^{(m)})_{n\in\N}]$ of sequences of integers where the partial quotients are computed by the Jacobi-Perron algorithm, obtained by iterating the following equations:
\begin{equation}\label{JPA}
\begin{cases}
& a_n^{(i)}=\lfloor x_n^{(i)}\rfloor ,\\
& x_{n+1}^{(1)}=\frac{1}{x_n^{(m)}-a_n^{(m)}},\\
& x_{n+1}^{(i)}=\frac{x_n^{(i-1)}-a_n^{(i-1)}}{x_n^{(m)}-a_n^{(m)}},\quad 2\leq i\leq m.
\end{cases}
\end{equation}
From \eqref{JPA}, it is possible to obtain the following identities for the complete quotients:
\begin{equation}\label{JPA2}
  x_n^{(i-1)} = a_n^{(i-1)} + \cfrac{x_{n+1}^{(i)}}{x_{n+1}^{(1)}}, \quad x_n^{(m)} = a_n^{(m)} + \cfrac{1}{x_{n+1}^{(1)}},  
\end{equation}
for $i=2, \ldots, m$ and $n = 0, 1 \ldots$, which provide the expansions of $(x^{(1)}_0,...,x^{(m)}_0)$ in a way analogue to \eqref{cf}.
Similarly to classical continued fractions, the $n$--th convergents of an MCF are given by the $m$-tuple $\left(\frac{A_n^{(1)}}{A_n^{(m+1)}},\ldots ,\frac{A_n^{(m)}}{A_n^{(m+1)}}\right)$, where
\begin{align*}
A_0^{(i)} & =a_0^{(i)},\quad A_0^{(m+1)}=1,\quad 1\leq i\leq m,\\
A_{-j}^{(i)} & =\delta_{i,j},\quad 1\leq i,j\leq m+1,\\
A_n^{(i)} & =\sum_{j=1}^m a_n^{(j)}A_{n-j}^{(i)}+A_{n-m-1}^{(i)},\quad 1\leq i\leq m+1.
\end{align*}
In the matrix form we may write
\begin{align}\label{convmat}
[A_{n-j}^{(i)}]_{1\leq i,j\leq m+1}=\prod_{k=0}^{n-1}\left[\begin{array}{cccc}
a_k^{(1)} & 1 & & \\
\vdots & & \ddots & \\
a_k^{(m)} & & & 1\\
1 & 0 & \cdots & 0
\end{array}\right].
\end{align}
If the Jacobi-Perron algorithm never stops (which holds, e.g., when the numbers $1,x^{(1)}_0,...,x^{(m)}_0$ are linearly independent over $\Q$), then
$$\lim_{n\to\infty}\frac{A_n^{(i)}}{A_n^{(m+1)}}=x^{(i)}_0$$
for each $i\in\{1,\ldots ,m\}$ (see, e.g., \cite{Dubois}).

The multidimensional continued fractions have been extensively studied from many points of view. First of all, there are several works involving the study of periodicity, since the Jacobi-Perron algorithm (and its variants) does not solve the Hermite problem, which is still an intriguing open problem in number theory. However, it has been proved the periodicity of the Jacobi-Perron algorithm for many families of algebraic irrationals, see, e.g., \cite{Ber1, Ber2, Ber3, Voutier, Rada, Lev}. Moreover, there are several works related to multidimensional continued fractions in Diophantine approximation for the study of simultaneous rational approximations of real numbers. A very good survey about this topic is \cite{Berthe}, where the authors discuss convergence properties, quality of approximations and dynamical properties. 

In this context, a study about the arithmetic of multidimensional continued fractions is still missing.
The aim of the paper is to present the generalization of Gosper's algorithm for multidimensional continued fractions. The paper is structured as follows. Section \ref{sec:uniGosper} presents a revisitation of classical Gosper's algorithm. In Section \ref{sec:MobiusMCF}, we study the M\"obius transforms of multidimensional continued fractions. In particular, we provide Algorithm \ref{Alg: GospMCF1} in order to compute, starting from the MCF of $(x^{(1)},,\ldots,x^{(m)})\in\R^m$, the MCF of
\begin{equation}\label{Eq: intro1}
\left(\frac{L^{(1)}({\bf x},1)}{L^{(m+1)}({\bf x},1)},\ldots ,\frac{L^{(m)}({\bf x},1)}{L^{(m+1)}({\bf x},1)}\right),
\end{equation}
where $L^{(1)},\ldots ,L^{(m+1)}:\R^{m+1}\to\R$ are linear forms.
In Section \ref{sec:arithMCF}, we define Algorithm \ref{Alg: GospMCF2} to compute the bilinear transforms of two MCFs, which allow to perform, in particular, sums and products of MCFs. Given $(x^{(1)},,\ldots,x^{(m)}),(y^{(1)},,\ldots,y^{(m)})\in\R^m$, Algorithm \ref{Alg: GospMCF2} provides the MCF of
\begin{equation}\label{Eq: intro2}
\left(\frac{L^{(1)}(({\bf x},1)({\bf y},1))}{L^{(m+1)}(({\bf x},1)({\bf y},1))},\ldots ,\frac{L^{(m)}(({\bf x},1)({\bf y},1))}{L^{(m+1)}(({\bf x},1)({\bf y},1))}\right),
\end{equation}
where $L^{(1)},\ldots ,L^{(m+1)}:\R^{m+1}\times\R^{m+1}\to\R$ are bilinear forms. In Section \ref{feas1} and Section \ref{feas2} we prove the correctness and feasibility of, respectively, Algorithm \ref{Alg: GospMCF1} and Algorithm \ref{Alg: GospMCF2}, i.e. that they end in a finite number of steps, providing the correct answer.\\
In Section \ref{sec:exp}, we present some numerical experiments for studying the behavior of the algorithm, mainly aimed to understand the number of inputs required to perform an output, i.e. how many partial quotients of the starting MCFs are required in order to compute correctly the partial quotients of the MCFs of \eqref{Eq: intro1} and \eqref{Eq: intro2}. The result is that the required number of inputs seems to grow always linearly, both for Algorithm \ref{Alg: GospMCF1} and Algorithm \ref{Alg: GospMCF2}, for any dimension $m$ of the MCF. Moreover, in all the cases they distribute around a line that has slope always less than $1$, suggesting that we require, on average, less than $1$ input in order to get $1$ output (see Section \ref{sec:exp} for more details). For future research, it would be interesting to prove that the number of inputs actually grows linearly with the number of outputs, and to obtain estimations about the slope of this line for the different cases. The SageMath code used for the analysis in Section \ref{sec:exp} is publicly available\footnote{\href{https://github.com/giulianoromeont/Continued-fractions}{https://github.com/giulianoromeont/Continued-fractions}}. Finally, in Section \ref{sec:mod} we present a possible improvement of the algorithm.

\section{Classical Gosper's algorithm revised}\label{sec:uniGosper}

First of all, given $x = [a_0; a_1, \ldots]$, Gosper's algorithm provides a method for obtaining the continued fraction expansion of $y = \frac{ax + b}{cx + d}$, with $a,b,c,d \in \mathbb Z$ and $ad - bc \not= 0$, so that we are able, e.g., to multiply and add rational numbers to a continued fraction.

If $\lfloor \frac{a}{c} \rfloor = \lfloor \frac{a+b}{c+d} \rfloor = b_0$, $x\geq 1$, and $c$, $c+d$ have the same sign, then we can easily see that $\lfloor y \rfloor = b_0$. Indeed, as $c$, $c+d$ have the same sign, we know that the function $cx+d$ has no zeros on the interval $[1,+\infty)$. Consequently, the function $\frac{ax+b}{cx+d}$ is well-defined, continuous and monotone on the interval $[1,+\infty)$. Moreover, the two values $\frac{a}{c}$ and $\frac{a+b}{c+d}$ represent the limit cases where, respectively, $x\rightarrow \infty$ and $x=1$. Summing up, if $x\in [1,+\infty)$, then $\frac{ax+b}{cx+d}$ lies between $\frac{a}{c}$ and $\frac{a+b}{c+d}$. Therefore, if $\lfloor \frac{a}{c} \rfloor = \lfloor \frac{a+b}{c+d} \rfloor$, we have no ambiguity in the choice of  $\lfloor \frac{ax+b}{cx+d} \rfloor$.\bigskip

According to the above, the main idea of the Gosper's algorithm is to exploit the partial quotients of $x$ in order to obtain such a situation where the first partial quotient of the continued fraction expansion of $y$ is known, i.e., $y = b_0 + \frac{1}{y_1}$. \smallskip 

We can use the partial quotients of $x$ observing that
\[ y = \cfrac{ax+b}{cx+d} = \cfrac{a\left(a_0+\frac{1}{x_1}\right)+b}{c\left(a_0+\frac{1}{x_1}\right)+d}=\cfrac{(aa_0+b)x_1 + a}{(ca_0+d)x_1 + c}, \]
where $x_1 = [a_1; a_2, \ldots]$, and in matrix form this is equivalent to:
\[\left[\begin{array}{cc}
a & b\\
c & d
\end{array}\right]\left[\begin{array}{cc}
a_0 & 1\\
1 & 0
\end{array}\right] = \left[\begin{array}{cc}
aa_0+b & a\\
ca_0+d & c
\end{array}\right].\]
Similarly, when we know that $y = b_0 + \frac{1}{y_1}$, we get
\[ y_1 = \cfrac{1}{y - b_0} = \cfrac{cx + d}{(a-cb_0)x + b -db_0} \]
and we can then focus on finding the floor of $y_1$ for obtaining the second partial quotients of the continued fraction expansion of $y$. Using matrices, the previous identity can also be written as
\[ \begin{bmatrix} 0 & 1 \cr 1 & - b_0 \end{bmatrix} \begin{bmatrix} a & b \cr c & d \end{bmatrix} = \begin{bmatrix} c & d \cr a - cb_0 & b- db_0 \end{bmatrix}. \]

In order to write a pseudo code of the Gosper's algorithm, it is useful to treat $x$ and $y$ as formal symbols that transform as follows when a partial quotient $p$ of $x$ is used (from now on by the symbol "$\mapsto$" we mean "is replaced by"):
\[x\mapsto p+\frac{1}{x}, \quad y=\frac{ax+b}{cx+d}\mapsto\frac{(ap+b)x+a}{(cp+d)x+c}\]
We may write it in the matrix form as
$$\left[\begin{array}{cc}
a & b\\
c & d
\end{array}\right]\mapsto\left[\begin{array}{cc}
ap+b & a\\
cp+d & c
\end{array}\right]=\left[\begin{array}{cc}
a & b\\
c & d
\end{array}\right]\left[\begin{array}{cc}
p & 1\\
1 & 0
\end{array}\right].$$
Similarly, when a partial quotient $q$ of $y$ is obtained, we have
$$y = \frac{ax+b}{cx+d}\mapsto\frac{cx+d}{(a-cq)x+(b-dq)}$$
and
$$\left[\begin{array}{cc}
a & b\\
c & d
\end{array}\right]\mapsto\left[\begin{array}{cc}
c & d\\
a-cq & b-dq
\end{array}\right]=\left[\begin{array}{cc}
0 & 1\\
1 & -q
\end{array}\right]\left[\begin{array}{cc}
a & b\\
c & d
\end{array}\right].$$
in the matrix form.

Then, Gosper's algorithm for computing the continued fraction expansion of $\frac{ax+b}{cx+d}$ is described in Algorithm \ref{gosper1}.

\IncMargin{1.5em}
\begin{algorithm}[h]\label{Alg: Gosp1}
	\caption{Gosper's algorithm for $\frac{ax+b}{cx+d}$}\label{gosper1}
	\SetKwData{Left}{left}
	\SetKwData{This}{this}
	\SetKw{And}{and}
	\SetKwFunction{Union}{Union}
	\SetKwFunction{FindCompress}{FindCompress}
	\SetKwInOut{Input}{Input}
	\SetKwInOut{Output}{Output}
    \SetKw{KwTo}{to}
    \SetKwComment{Comment}{$\triangleright$ }{ }
	\Input{$x = [a_0; a_1, \ldots]$, $a,b,c,d \in \mathbb Z$, $M,N \in \mathbb N$}
	\Output{$[b_0; b_1, \ldots] = \frac{ax+b}{cx+d}$}
	\BlankLine
	$r \gets 0$, $s \gets 0$, $t \gets 0$ 

    \While{$s<M$ \And $r<N$}{

        \If{$\sgn (c)=\sgn (c+d)$ and $\lfloor \frac{a}{c} \rfloor = \lfloor \frac{a+b}{c+d} \rfloor$}{
            $b_s \gets \lfloor \frac{a}{c} \rfloor$
        
            $\begin{bmatrix} a & b \cr c & d \end{bmatrix} \gets \begin{bmatrix} 0 & 1 \cr 1 & -b_s \end{bmatrix}\begin{bmatrix} a & b \cr c & d \end{bmatrix}$

            $s \gets s + 1$
        }

        \Else{
            $\begin{bmatrix} a & b \cr c & d \end{bmatrix} \gets \begin{bmatrix} a & b \cr c & d \end{bmatrix}\begin{bmatrix} a_t & 1 \cr 1 & 0 \end{bmatrix}$

            $t \gets t + 1$
        }

    $r \gets r+1$

    }
    
    \Return $[b_0; b_1, \ldots , b_{s-1}]$

\end{algorithm}
\DecMargin{1.5em}





In Algorithm \ref{Alg: Gosp1}, the index $r$, that is bounded by $N$, counts the total number of times that the algorithm is executed, either performing input or output transformations. The reason to bound $r$ is that, a priori, we do not know how many input transformations we need before providing an output. Therefore, if we gave a bound only on the desired number of output partial quotients, the algorithm could run for an arbitrarily large number of steps.

\begin{Example}
Let us consider $\sqrt{2}=[1,\overline{2}]$ and let us suppose that we want to compute the continued fraction expansion of $2\sqrt{2}$. In this case
\begin{equation}\label{Eq: 2sqrt2}
    2\sqrt{2}=\frac{2\sqrt{2}+0}{0+1},
\end{equation}
so that $a=2,b=c=0,d=1$. The first partial quotient is $q=\lfloor \sqrt{2}\rfloor=1$, and we substitute $\sqrt{2}=1+\frac{1}{x_1}$ in \eqref{Eq: 2sqrt2}, in order to obtain:
\begin{equation*}
    2\sqrt{2}=\frac{2\left(1+\frac{1}{x_1}\right)+0}{0+1}=\frac{2x_1+2}{x_1}.
\end{equation*}
Since
\[\left\lfloor\frac{2}{1}\right\rfloor=2\neq 4=\left\lfloor\frac{2+2}{1}\right\rfloor,\]
we can not decide yet which is the partial quotient. The second partial quotient is $\lfloor x_1\rfloor=2$, so that
\begin{equation*}
    2\sqrt{2}=\frac{2\left(2+\frac{1}{x_2}\right)+2}{\left(2+\frac{1}{x_2}\right)}=\frac{6x_2+2}{2x_2+1}.
\end{equation*}
Since
\[\left\lfloor\frac{6}{2}\right\rfloor=3\neq2=\left\lfloor\frac{6+2}{2+1}\right\rfloor,\]
we continue, i.e.
\begin{equation*}
2\sqrt{2}=\frac{6\left(2+\frac{1}{x_3}\right)+2}{2\left(2+\frac{1}{x_3}\right)+1}=\frac{14x_3+6}{5x_3+2}.
\end{equation*}
Now we have that
\[\left\lfloor\frac{14}{5}\right\rfloor=2=\left\lfloor\frac{14+6}{5+2}\right\rfloor,\]
and the first partial quotient is then $2$. Using the matrix notation of Algorithm \ref{Alg: Gosp1}, we have computed the product
\[\begin{bmatrix} 2 & 0 \cr 0 & 1 \end{bmatrix}\begin{bmatrix} 1 & 1 \cr 1 & 0 \end{bmatrix}\begin{bmatrix} 2 & 1 \cr 1 & 0 \end{bmatrix}\begin{bmatrix} 2 & 1 \cr 1 & 0 \end{bmatrix}=\begin{bmatrix} 14 & 6 \cr 5 & 2 \end{bmatrix},\]
until an output is reached, i.e. $b_0=\lfloor 2\sqrt{2} \rfloor=2$.
Therefore we now have
\[\frac{1}{2\sqrt{2}-2},\]
i.e. a new transformation where $a=0,b=1,c=2,d=-2$. In matrix form it is
\[\begin{bmatrix} 0 & 1 \cr 1 & -2 \end{bmatrix}\begin{bmatrix} 2 & 0 \cr 0 & 1 \end{bmatrix}=\begin{bmatrix} 0 & 1 \cr 2 & -2 \end{bmatrix}.\]
However, we already computed the first $3$ input matrices, so that we avoid to start again from the first partial quotient of $\sqrt{2}$, computing directly
\[\begin{bmatrix} 0 & 1 \cr 1 & -2 \end{bmatrix}\begin{bmatrix} 14 & 6 \cr 5 & 2 \end{bmatrix}=\begin{bmatrix} 5 & 2 \cr 4 & 2 \end{bmatrix}.\]
Then we keep performing the algorithm on the latter matrix.
\end{Example}

Now, for computing the continued fraction expansion of the sum, difference, product or quotient of two continued fractions $x= [a_0; a_1, \ldots]$ and $y = [b_0; b_1, \ldots]$, we need to deal with their bilinear fractional transformation
\begin{align}\label{genexpr}
    z = \frac{axy+bx+cy+d}{exy+fx+gy+h},\quad a,b,c,d,e,f,g,h\in\Z .
\end{align}
Similarly to the case of linear fractional transformations, if
\[\left\lfloor \frac{a}{e}\right\rfloor = \left\lfloor \frac{a+b}{e+f}\right\rfloor = \left\lfloor \frac{a+c}{e+g}\right\rfloor = \left\lfloor \frac{a+b+c+d}{e+f+g+h}\right\rfloor = c_0,\]
and $e$, $e+f$, $e+g$, $e+f+g+h$ have the same sign, then we have $\lfloor z \rfloor = c_0$ and $z = c_0 + \frac{1}{z_1}$. This follows from the fact that then the denominator $exy+fx+gy+h$ has no zeros for $x,y\in [1,+\infty)$ and then $\frac{axy+bx+cy+d}{exy+fx+gy+h}$ is a well-defined continuous function on $[1,+\infty)^2$ and monotone with respect to each variable. Thus, we can exploit the partial quotients of $x$ and $y$ in order to get to such situation. In particular, we have
\[ z = \cfrac{(aa_0+c)x_1y + (ba_0+d)x_1 + ay + b}{(ea_0+g)x_1y + (fa_0+h)x_1 + ey + f}, \quad z = \cfrac{(ab_0+b)xy_1 + ax + (cb_0+d)y_1 + c}{(eb_0+f)xy_1 +ex + (gb_0+h)y_1+g}, \]
where $x_1 = [a_1; a_2, \ldots]$ and $y_1 = [b_1; b_2, \ldots]$. On the other hand, if $z = c_0 + \frac{1}{z_1}$, we obtain
\[ z_1 = \cfrac{1}{z - c_0} = \cfrac{exy + fx + gy + h}{(a-ec_0)xy + (b-fc_0)x + (c-gc_0)y + d - hc_0}. \]
We summarize this procedure in Algorithm \ref{Alg: Gosp2}.

\IncMargin{1.5em}
\begin{algorithm}[h]\label{Alg: Gosp2}
	\caption{Gosper's algorithm for $\frac{axy+bx+cy+d}{exy+fx+gy+h}$}\label{gosper1}
	\SetKwData{Left}{left}
	\SetKwData{This}{this}
	\SetKw{And}{and}
	\SetKwFunction{Union}{Union}
	\SetKwFunction{FindCompress}{FindCompress}
	\SetKwInOut{Input}{Input}
	\SetKwInOut{Output}{Output}
    \SetKw{KwTo}{to}
    \SetKwComment{Comment}{$\triangleright$ }{ }
	\Input{$x = [a_0; a_1, \ldots]$, $y = [b_0; b_1, \ldots]$, $a,b,c,d,e,f,g,h \in \mathbb Z$, $M,N \in \mathbb N$}
	\Output{$[c_0; c_1, \ldots] = \frac{axy+bx+cy+d}{exy+fx+gy+h}$}
	\BlankLine
	$r \gets 0$, $s \gets 0$, $t \gets 0$, $u \gets 0$, $v \gets 0$

    \While{$s<M$ \And $r<N$}{

        \If{$\sgn(e)=\sgn(e+f)=\sgn(e+g)=\sgn(e+f+g+h)$ and $\lfloor \frac{a}{e}\rfloor = \lfloor \frac{a+b}{e+f}\rfloor = \lfloor \frac{a+c}{e+g}\rfloor = \lfloor \frac{a+b+c+d}{e+f+g+h}\rfloor$}{
            $c_s \gets \lfloor \frac{a}{e} \rfloor$
        
            $\begin{bmatrix} a & b & c & d \cr e & f & g & h \end{bmatrix} \gets \begin{bmatrix} e & f & g & h \cr a-ec_s & b-fc_s & c - gc_s & d-hc_s \end{bmatrix}$

            $s \gets s + 1$
        }

        \Else{
        
        \If{$u = v$}{
              $\begin{bmatrix} a & b & c & d \cr e & f & g & h \end{bmatrix} \gets \begin{bmatrix} aa_u+c & ba_u+d & a & b \cr ea_u + g & fa_u+h & e & f \end{bmatrix}$ 
              
              $u \gets u + 1$
              } 
              
            \Else{$\begin{bmatrix} a & b & c & d \cr e & f & g & h \end{bmatrix} \gets \begin{bmatrix} ab_v+b & a & cb_v+d & c \cr eb_v+f & e & gb_v+h & g \end{bmatrix}$

            $v \gets v+1$
            }

            $t \gets t+1$
        }
            
        $r \gets r+1$

    }
    
    \Return $[c_0; c_1, \ldots ,c_{s-1}]$

\end{algorithm}
\DecMargin{1.5em}

\section{M\"obius transform of multidimensional continued fractions.} \label{sec:MobiusMCF}

In this section, we aim to define the analogue of Algorithm \ref{Alg: Gosp1} for computing the M\"{o}bius transformation of a single MCF. First, let us focus on a $m$-tuple ${\bf x}=(x^{(1)},\ldots ,x^{(m)})$ whose MCF expansion is given by $$\left[\textbf{a}^{(i)}\right]_{i=1}^m=[(a_n^{(1)})_{n\in\N},...,(a_n^{(m)})_{n\in\N}].$$ 

In the one-dimensional case, the Gosper's algorithm provides the continued fraction expansion of $\frac{ax+b}{cx+d}$, $a,b,c,d\in\Z$, where $\frac{ax+b}{cx+d}$ is a M\"{o}bius transformation of $x$ with integral coefficients. Hence, if we identify $x$ with a point $[x:1]$ in the projective line $\mathbb{P}^1(\R)$, then $\frac{ax+b}{cx+d}$ identifies with $[ax+b:cx+d]\in\mathbb{P}^1(\R)$, which is a projective transformation of $[x:1]$. This can be easily extended into higher dimensions. Namely, identify ${\bf x}$ with $[x^{(1)}:\ldots :x^{(m)}:1]\in\mathbb{P}^m(\R)$. Its projective transformation in $\mathbb{P}^m(\R)$ has the form $$[L^{(1)}({\bf x},1):\ldots :L^{(m)}({\bf x},1):L^{(m+1)}({\bf x},1)],$$ where $L^{(1)},\ldots ,L^{(m+1)}:\R^{m+1}\to\R$ are linear forms.\\
The point $[L^{(1)}({\bf x},1):\ldots :L^{(m)}({\bf x},1):L^{(m+1)}({\bf x},1)]\in\mathbb{P}^m(\R)$ corresponds to the point $$\left(\frac{L^{(1)}({\bf x},1)}{L^{(m+1)}({\bf x},1)},\ldots ,\frac{L^{(m)}({\bf x},1)}{L^{(m+1)}({\bf x},1)}\right)$$ in the affine space $\R^m$. This point, where $L^{(1)},\ldots ,L^{(m+1)}$ have all integral coefficients, is our candidate for an $m$-tuple, for which we will compute its MCF expansion.

In the higher-dimensional case, if the Jacobi-Perron algorithm stops at the $n$-th step, i.e. $x_n^{(m)}\in\Z$ and $x_k^{(m)}\not\in\Z$ for $k\in\{0,1,\ldots ,n\}$, it is possible that $\frac{A_{n+1}^{(i)}}{A_{n+1}^{(m+1)}}\neq x^{(i)}$ for some $i\in\{1,\ldots ,m\}$, see, e.g., the MCF expansion of $(\pi,5\pi+2)$. Thus, we restrict ourselves to the case when both $(x^{(1)},\ldots x^{(m)},1)$ and $\left(\frac{L^{(1)}({\bf x},1)}{L^{(m+1)}({\bf x},1)},\ldots ,\frac{L^{(m)}({\bf x},1)}{L^{(m+1)}({\bf x},1)},1\right)$ are linearly independent over $\Q$. This guarantees that the Jacobi-Perron algorithm for $(x^{(1)},\ldots ,x^{(m)})$ and $\left(\frac{L^{(1)}({\bf x},1)}{L^{(m+1)}({\bf x},1)},\ldots ,\frac{L^{(m)}({\bf x},1)}{L^{(m+1)}({\bf x},1)}\right)$ provides convergents that converge to them. Note that $\Q$-linear independence of the $m+1$-tuple $\left(\frac{L^{(1)}({\bf x},1)}{L^{(m+1)}({\bf x},1)},\ldots ,\frac{L^{(m)}({\bf x},1)}{L^{(m+1)}({\bf x},1)},1\right)$ is equivalent to $\Q$-linear independence of the $m+1$-tuple $(L^{(1)}({\bf x},1),\ldots ,L^{(m)}({\bf x},1),L^{(m+1)}({\bf x},1))$. This, with the assumption of $\Q$-linear independence of the $m+1$-tuple $(x_1,\ldots x_m,1)$, is equivalent to $\Q$-linear independence of the linear forms $L^{(1)},\ldots ,L^{(m)},L^{(m+1)}$.

Summing up, our input is the continued fraction expansion $\textbf{a}=\left(\textbf{a}^{(i)}\right)_{i=1}^m=\left(\left(a_k^{(i)}\right)_{k=0}^\infty\right)_{i=1}^m$ of ${\bf x}=\left(x^{(i)}\right)_{i=1}^m$ such that $(x^{(1)},\ldots x^{(m)},1)$ is linearly independent over $\Q$, and
\[C=\left[c_{j}^{(i)}\right]_{1\leq i,j\leq m+1}\in M(m+1\times m+1,\Z)\]
such that $\det C\neq 0$. Our output is the continued fraction expansion $\left(\textbf{b}^{(i)}\right)_{i=1}^m=\left(\left(b_k^{(i)}\right)_{k=0}^\infty\right)_{i=1}^m$ of $\left(\frac{L_1({\bf x},1)}{L_{m+1}({\bf x},1)},\ldots ,\frac{L_m({\bf x},1)}{L_{m+1}({\bf x},1)}\right)$, where
$$L^{(i)}({\bf x},1)=\sum_{j=1}^m c_j^{(i)}x_j+c_{m+1}^{(i)},\ i\in\{1,\ldots m+1\}.$$

Remembering that
\begin{equation*}
  \begin{cases}
    x_k^{(m)} & =a_k^{(m)}+\frac{1}{x_{k+1}^{(1)}},\\
    x_k^{(i-1)} & =a_k^{(i-1)}+\frac{x_{k+1}^{(i)}}{x_{k+1}^{(1)}},\quad 2\leq i\leq m
  \end{cases}  
\end{equation*}
we could think of ${\bf x}$ as an $m$-tuple of formal symbols that transform in the following way when we exploit the first partial quotient:
$$x^{(m)}\mapsto a_0^{(m)}+\frac{1}{x^{(1)}},\quad x^{(i)}\mapsto a_0^{(i)}+\frac{x^{(i+1)}}{x^{(1)}},\quad i\in\{1,\ldots m-1\}.$$
Then, we also have
\begin{equation}\label{intrans}
  \begin{split}
\frac{L^{(i)}({\bf x},1)}{L^{(m+1)}({\bf x},1)}\mapsto
&\frac{L^{(i)}\left(a_0^{(1)}+\frac{x^{(2)}}{x^{(1)}},\ldots ,a_0^{(m-1)}+\frac{x^{(m)}}{x^{(1)}},a_0^{(m)}+\frac{1}{x^{(1)}},1\right)}{L^{(m+1)}\left(a_0^{(1)}+\frac{x^{(2)}}{x^{(1)}},\ldots ,a_0^{(m-1)}+\frac{x^{(m)}}{x^{(1)}},a_0^{(m)}+\frac{1}{x^{(1)}},1\right)}=\\
=&\frac{L^{(i)}(a_0^{(1)}x^{(1)}+x^{(2)},\ldots ,a_0^{(m-1)}x^{(1)}+x^{(m)},a_0^{(m)}x^{(1)}+1,x^{(1)})}{L^{(m+1)}(a_0^{(1)}x^{(1)}+x^{(2)},\ldots ,a_0^{(m-1)}x^{(1)}+x^{(m)},a_0^{(m)}x^{(1)}+1,x^{(1)})},\ i\in\{1,\ldots m\}.
  \end{split}  
\end{equation}
If we write 
$$L^{(i)}(x^{(1)},\ldots ,x^{(m+1)})=\sum_{j=1}^{m+1} c_j^{(i)}x^{(j)},\ i\in\{1,\ldots m+1\},$$
and $C=\left[c_{j}^{(i)}\right]_{1\leq i,j\leq m+1}$, then by the transformation \eqref{intrans} we get the following transformation of the matrix $C$:
\begin{equation}\label{matintrans}
C\mapsto C\left[\begin{array}{cccccc}
a_0^{(1)} & 1 & \cdots & 0 & 0 \\
\vdots & \vdots & \ddots & \vdots & \vdots\\
a_0^{(m-1)} & 0 & \cdots & 1 & 0\\
a_0^{(m)} & 0 & \cdots & 0 & 1\\
1 & 0 & \cdots & 0 & 0
\end{array}\right].
\end{equation}

Similarly, we can see how the M\"obius transform of an MCF changes when we know the first partial quotient of its expansion. 
For each $i\in\{1,\ldots ,m\}$ the value of $\frac{L^{(i)}({\bf x},1)}{L^{(m+1)}({\bf x},1)}$ lies between the minimum and maximum of the values $\frac{c_j^{(i)}}{c_j^{(m+1)}}$, $j\in\{1,\ldots ,m+1\}$, on condition that all the values $c_j^{(m+1)}$, $j\in\{1,\ldots ,m+1\}$, have the same sign. This is true in virtue of the following mediant inequality
\begin{equation}\label{med}
    \frac{a}{b}\leq\frac{a+c}{b+d}\leq\frac{c}{d},\, \text{ when } \frac{a}{b}\leq\frac{c}{d} \text{ and } bd>0 
\end{equation}
and induction with respect to $m$. Thus, if all the mentioned values have the same integral part, then also $\frac{L^{(i)}({\bf x},1)}{L^{(m+1)}({\bf x},1)}$ has this value as the integral part. Then, we have
\[b^{(i)}_0=\left\lfloor\frac{L^{(i)}({\bf x},1)}{L^{(m+1)}({\bf x},1)}\right\rfloor=\left\lfloor\frac{c_1^{(i)}}{c_1^{(m+1)}}\right\rfloor.\]
If, for each $i\in\{1,\ldots m\}$ and $j_1,j_2\in\{1,\ldots m+1\}$,
\begin{equation}\label{Eq: floors}
\left\lfloor\frac{c_{j_1}^{(i)}}{c_{j_1}^{(m+1)}}\right\rfloor =\left\lfloor\frac{c_{j_2}^{(i)}}{c_{j_2}^{(m+1)}}\right\rfloor
\end{equation}
then we know the first partial quotient of the MCF expansion of $\left(\frac{L_1({\bf x},1)}{L_{m+1}({\bf x},1)},\ldots ,\frac{L_m({\bf x},1)}{L_{m+1}({\bf x},1)}\right)$. Notice that, if \eqref{Eq: floors} is not satisfied for all $i$, we must keep getting inputs from all the $x^{(i)} $, as all $x^{(i)} $ appear in all linear forms.  
Then, we can write
\begin{align*}
 \frac{L^{(1)}(x,y)}{L^{(m+1)}(x,y)} & \mapsto\frac{1}{\frac{L^{(m)}(x,y)}{L^{(m+1)}(x,y)}-b^{(m)}}=\frac{L^{(m+1)}(x,y)}{L^{(m)}(x,y)-b^{(m)}L^{(m+1)}(x,y)},\\
\frac{L^{(i+1)}(x,y)}{L^{(m+1)}(x,y)} & \mapsto\frac{\frac{L^{(i)}(x,y)}{L^{(m+1)}(x,y)}-b^{(i)}}{\frac{L^{(m)}(x,y)}{L^{(m+1)}(x,y)}-b^{(m)}}=\frac{L^{(i)}(x,y)-b^{(i)}L^{(m+1)}(x,y)}{L^{(m)}(x,y)-b^{(m)}L^{(m+1)}(x,y)},\, i\in\{1,\ldots m-1\}.
\end{align*}
This means that the matrix $C$ transforms as follows, after getting an output:
$$C\mapsto\left[\begin{array}{ccccc}
0 & 0 & \cdots & 0 & 1\\
1 & 0 & \cdots & 0 & -b^{(1)}\\
0 & 1 & \cdots & 0 & -b^{(2)}\\
\vdots & \vdots & \ddots & \vdots & \vdots\\
0 & 0 & \cdots & 1 & -b^{(m)}
\end{array}\right] C.$$

We summarize this procedure in Algorithm \ref{Alg: GospMCF1}, and it is then illustrated in Example \ref{E1} in the Appendix.

\IncMargin{1.5em}
\begin{algorithm}[h]\label{Alg: GospMCF1}
	\caption{Algorithm for obtaining the MCF expansion of $\left(\frac{L_1({\bf x},1)}{L_{m+1}({\bf x},1)},\ldots ,\frac{L_m({\bf x},1)}{L_{m+1}({\bf x},1)}\right)$ given the MCF expansion of $\mathbf x = (x^{(1)}, \ldots, x^{(m)})$, where $L_1, \ldots, L_{m+1}$ are linear forms.}\label{gosper1}
	\SetKwData{Left}{left}
	\SetKwData{This}{this}
	\SetKw{And}{and}
	\SetKwFunction{Union}{Union}
	\SetKwFunction{FindCompress}{FindCompress}
	\SetKwInOut{Input}{Input}
	\SetKwInOut{Output}{Output}
    \SetKw{KwTo}{to}
    \SetKwComment{Comment}{$\triangleright$ }{ }
	\Input{$\textbf{x}= [(a_n^{(1)})_{n\geq0},...,(a_n^{(m)})_{n\geq0}]$, $C=\left[c_{j}^{(i)}\right]_{1\leq i,j\leq m+1}$, $M,N \in \mathbb N$}
	\Output{$[(b_n^{(1)})_{n\geq0},...,(b_n^{(m)})_{n\geq0}] = \left(\frac{L_1({\bf x},1)}{L_{m+1}({\bf x},1)},\ldots ,\frac{L_m({\bf x},1)}{L_{m+1}({\bf x},1)}\right)$}
	\BlankLine
	$r \gets 0$, $s \gets 0$, $t \gets 0$

    \While{$s<M$ \And $r<N$}{

        \If{$\sgn(c_{j_1}^{(m+1)})=\sgn(c_{j_2}^{(m+1)})$ and $\left\lfloor\frac{c_{j_1}^{(i)}}{c_{j_1}^{(m+1)}}\right\rfloor =\left\lfloor\frac{c_{j_2}^{(i)}}{c_{j_2}^{(m+1)}}\right\rfloor$ for each $i\in\{1,\ldots m\}$ and $j_1,j_2\in\{1,\ldots m+1\}$}{
            $b^{(i)}_s \gets \left\lfloor\frac{c_{1}^{(i)}}{c_{1}^{(m+1)}}\right\rfloor ,\, i\in\{1,\ldots m\},$
            
            $C\gets\left[\begin{array}{ccccc}
            0 & 0 & \cdots & 0 & 1\\
            1 & 0 & \cdots & 0 & -b^{(1)}_s\\
            0 & 1 & \cdots & 0 & -b^{(2)}_s\\
            \vdots & \vdots & \ddots & \vdots & \vdots\\
            0 & 0 & \cdots & 1 & -b^{(m)}_s
            \end{array}\right] C,$
        
            $s \gets s + 1$
        }

        \Else{
            $C\gets C\left[\begin{array}{cccccc}
            a_t^{(1)} & 1 & \cdots & 0 & 0 \\
            a_t^{(2)} & 0 & 1 & \cdots & 0 \\
            \vdots & \vdots & \ddots & \vdots & \vdots\\
            a_t^{(m-1)} & 0 & \cdots & 1 & 0\\
            a_t^{(m)} & 0 & \cdots & 0 & 1\\
            1 & 0 & \cdots & 0 & 0
            \end{array}\right]$

            $t \gets t + 1$
        }

    $r \gets r+1$

    }
    
    \Return $[(b_n^{(1)})_{n=0}^{s-1},...,(b_n^{(m)})_{n=0}^{s-1}]$

\end{algorithm}
\DecMargin{1.5em}

\newpage
\subsection{Proof of feasibility of the algorithm}\label{feas1}

We need to prove the validity of the multidimensional version of the Gosper's algorithm, i.e. that after using a finite number of partial quotients of the MCF of $\mathbf x$ we can obtain a partial quotient of the MCF of the M\"obius transform of $\mathbf x$.

We see that the $n$-fold iteration of the transformation \eqref{matintrans} of the matrix $C$ is
\begin{equation*}
C\mapsto C\cdot\prod_{k=0}^{n-1}\left[\begin{array}{cccccc}
a_k^{(1)} & 1 & \cdots & 0 & 0 \\
\vdots & \vdots & \ddots & \vdots & \vdots\\
a_k^{(m-1)} & 0 & \cdots & 1 & 0\\
a_k^{(m)} & 0 & \cdots & 0 & 1\\
1 & 0 & \cdots & 0 & 0
\end{array}\right]=C\cdot [A_{n-j}^{(i)}]_{1\leq i,j\leq m+1},
\end{equation*}
where the equality comes from \eqref{convmat}. This means that meanwhile the linear forms $L^{(i)}$, $i\in\{1,\ldots ,m+1\}$, transform as follows:
\begin{equation*}
  \begin{split}
L^{(i)}({\bf x},1)\mapsto&L^{(i)}\left(\sum_{j=1}^{m}A_{n+1-j}^{(1)}x^{(j)}+A_{n-m}^{(1)},\ldots ,\sum_{j=1}^{m}A_{n+1-j}^{(m+1)}x^{(j)}+A_{n-m}^{(m+1)}\right)\\
&=\sum_{j=1}^m\left(\sum_{k=0}^{m+1} c_k^{(i)}A_{n+1-j}^{(k)}\right)x^{(j)}+\sum_{k=0}^{m+1} c_k^{(i)}A_{n-m}^{(k)},\, i\in\{1,\ldots m\}.
  \end{split}  
\end{equation*}
Hence,
\begin{equation}\label{coefftrans}
    c_j^{(i)}\mapsto\sum_{k=0}^{m+1} c_k^{(i)}A_{n+1-j}^{(k)},\, i,j\in\{1,\ldots ,m+1\}
\end{equation}
and
\begin{equation}\label{ratiotrans}
  \begin{split}
\frac{L^{(i)}({\bf x},1)}{L^{(m+1)}({\bf x},1)}\mapsto\frac{\sum_{j=1}^m\left(\sum_{k=0}^{m+1} c_k^{(i)}A_{n+1-j}^{(k)}\right)x^{(j)}+\sum_{k=0}^{m+1} c_k^{(i)}A_{n-m}^{(k)}}{\sum_{j=1}^m\left(\sum_{k=0}^{m+1} c_k^{(m+1)}A_{n+1-j}^{(k)}\right)x^{(j)}+\sum_{k=0}^{m+1} c_k^{(m+1)}A_{n-m}^{(k)}},\, i\in\{1,\ldots m\}.
  \end{split}  
\end{equation}
As all the values of $x^{(j)}$ are positive, the last expression in \eqref{ratiotrans} lies between the minimum and maximum of the values 
\begin{equation}\label{coeffratiotrans}
    \frac{\sum_{k=0}^{m+1} c_k^{(i)}A_{n+1-j}^{(k)}}{\sum_{k=0}^{m+1} c_k^{(m+1)}A_{n+1-j}^{(k)}}=\frac{\sum_{k=0}^{m} c_k^{(i)}\frac{A_{n+1-j}^{(k)}}{A_{n+1-j}^{(m+1)}}+c_{m+1}^{(i)}}{\sum_{k=0}^{m} c_k^{(m+1)}\frac{A_{n+1-j}^{(k)}}{A_{n+1-j}^{(m+1)}}+c_{m+1}^{(m+1)}},\, j\in\{1,\ldots ,m+1\},
\end{equation}
on condition that all of the values $\sum_{k=0}^{m+1}c_k^{(m+1)}A_{n+1-j}^{(k)}$, $j\in\{1,\ldots ,m+1\}$, have the same sign. Since the $m+1$-tuple $(x_0^{(1)},\ldots ,x_0^{(m)},1)$ is $\Q$-linearly independent, we know that 
$$\lim_{n\to\infty}\frac{A_{n+1-j}^{(k)}}{A_{n+1-j}^{(m+1)}}=x_0^{(k)},$$ for $k\in\{1,\ldots m\}$, $j\in\{1,\ldots ,m+1\}$. From this, \eqref{coefftrans}, and \eqref{coeffratiotrans} we see that the consecutive transformations of the ratios $\frac{c_j^{(i)}}{c_j^{(m+1)}}$, $j\in\{1,\ldots ,m+1\}$, converge to the same value, namely $\frac{L^{(i)}(x_0^{(1)},\ldots ,x_0^{(m)},1)}{L^{(m+1)}(x_0^{(1)},\ldots ,x_0^{(m)},1)}$. In particular, for sufficiently large $n\in\N$ all the values 
$$\sum_{k=0}^{m} c_k^{(m+1)}\frac{A_{n+1-j}^{(k)}}{A_{n+1-j}^{(m+1)}}+c_{m+1}^{(m+1)} \quad j\in\{1,\ldots ,m+1\}$$
have the same sign. Consequently, all the values $\sum_{k=0}^{m+1}c_k^{(m+1)}A_{n+1-j}^{(k)}$, $j\in\{1,\ldots ,m+1\}$, have the same sign as all of the values $A_j^{(m+1}$, $j\in\N$, are positive. Since the $m+1$-tuple $\left(\frac{L^{(1)}(x_0^{(1)},\ldots ,x_0^{(m)},1)}{L^{(m+1)}(x_0^{(1)},\ldots ,x_0^{(m)},1)},\ldots ,\frac{L^{(m)}(x_0^{(1)},\ldots ,x_0^{(m)},1)}{L^{(m+1)}(x_0^{(1)},\ldots ,x_0^{(m)},1)},1\right)$ is $\Q$-linearly independent, for each\linebreak $i\in\{1,\ldots ,m\}$ the value $\frac{L^{(i)}(x_0^{(1)},\ldots ,x_0^{(m)},1)}{L^{(m+1)}(x_0^{(1)},\ldots ,x_0^{(m)},1)}$ is irrational, in particular non-integral. This means that for sufficiently large $n$, the $n$-th iterations of the ratios $\frac{c_j^{(i)}}{c_j^{(m+1)}}$, $j\in\{1,\ldots ,m+1\}$, have the same integral part and the $n$-th iterations of the coefficients $c_j^{(i)}$, $j\in\{1,\ldots ,m+1\}$, have the same sign. Hence, we may perform output transformation.

There remains an open problem on the number of input transformations sufficient to perform the output transformation. This problem is relevant from the computational complexity point of view. From the above reasoning, we should expect that this number depends essentially on the values $x_0^{(1)},\ldots ,x_0^{(m)}$ (in particular, the rate of convergence of their convergents) and the entries of the matrix $C$. In the case of $m=1$, i.e. classical continued fractions, Liardet and Stambul gave upper bounds for this number in terms of the matrix $C$ (see \cite[Lemma 6 and Theorem 2]{LiarStam}).

\section{Sums and products of MCFs}\label{sec:arithMCF}

Now we work with a pair of $m$-tuples ${\bf x}=(x^{(1)},\ldots ,x^{(m)})$ and ${\bf y}=(y^{(1)},\ldots ,y^{(m)})$ aiming at performing operations between two MCFs.
Once again we shall use linear algebra and projective geometry to motivate our choice of possible output. In the classical, one-dimensional case we have the continued fractions of $x,y\in\R$ and $a,b,c,d,e,f,g,h\in\Z$ as the input and the continued fraction of $\frac{axy+bx+cy+d}{exy+fx+gy+h}$ as the output. After identification with projective points we get $[x:1]$ and $[y:1]$ as a part of the input and $[{axy+bx+cy+d}:{exy+fx+gy+h}]$ as the output. Here we note that $({axy+bx+cy+d},{exy+fx+gy+h})$ is a bilinear transformation of the pair of vectors $(x,1)$, $(y,1)$. In higher dimensions this generalizes to $[x^{(1)}:\ldots :x^{(m)}:1]$ and $[y^{(1)}:\ldots :y^{(m)}:1]$ as a part of the input and $[L^{(1)}(({\bf x},1),({\bf y},1)):\ldots :L^{(m)}(({\bf x},1),({\bf y},1)):L^{(m+1)}(({\bf x},1),({\bf y},1))]$ as the output, where $L^{(1)},\ldots ,L^{(m+1)}$ are bilinear forms with integral coefficients (recall that $({\bf x},1)=(x^{(1)},\ldots ,x^{(m)},1)$ and $({\bf y},1)=(y^{(1)},\ldots ,y^{(m)},1)$). Coming back to affine space and continued fractions, we get the MCF expansions $\textbf{a}=\left(\textbf{a}^{(i)}\right)_{i=1}^m=\left(\left(a_k^{(i)}\right)_{k=0}^\infty\right)_{i=1}^m$ and $\textbf{b}=\left(\textbf{b}^{(i)}\right)_{i=1}^m=\left(\left(b_k^{(i)}\right)_{k=0}^\infty\right)_{i=1}^m$ of ${\bf x}=(x_1,\ldots ,x_m)$ and ${\bf y}=(y_1,\ldots ,y_m)\in\R^m$, respectively, and matrices $C^{(i)}=\left[c_{j,l}^{(i)}\right]_{1\leq j,l\leq m+1}\in M(m+1\times m+1,\Z)$, $i\in\{1,\ldots m+1\}$, as the input and the continued fraction expansion $\textbf{d}=\left(\textbf{d}^{(i)}\right)_{i=1}^m=\left(\left(d_k^{(i)}\right)_{k=0}^\infty\right)_{i=1}^m$ of $\left(\frac{L^{(1)}(({\bf x},1)({\bf y},1))}{L^{(m+1)}(({\bf x},1)({\bf y},1))},\ldots ,\frac{L^{(m)}(({\bf x},1)({\bf y},1))}{L^{(m+1)}(({\bf x},1)({\bf y},1))}\right)$, where 
$$L^{(i)}((x^{(1)},\ldots x^{(m+1)}),(y^{(1)},\ldots y^{(m+1)}))=\sum_{1\leq j,l\leq m+1} c_{j,l}^{(i)}x^{(j)}y^{(l)},\ i\in\{1,\ldots m+1\},$$
as the output.

Just as in the case of the single $m$-tuple ${\bf x}$ we additionally assume that all the mentioned continued fractions are infinite in order to have the convergence of the Jacobi-Perron algorithm. Moreover, we assume that the $m+1$-tuple $\left(\frac{L^{(1)}(({\bf x},1)({\bf y},1))}{L^{(m+1)}(({\bf x},1)({\bf y},1))},\ldots ,\frac{L^{(m)}(({\bf x},1)({\bf y},1))}{L^{(m+1)}(({\bf x},1)({\bf y},1))},1\right)$ is $\Q$-linearly independent.


When we exploit the first partial quotient of the MCF expansion of ${\bf x}$, it transforms as
$$x^{(m)}\mapsto a_0^{(m)}+\frac{1}{x^{(1)}}, x^{(i)}\mapsto a_0^{(i)}+\frac{x^{(i+1)}}{x^{(1)}}, i\in\{1,\ldots m-1\},$$
and consequently
\begin{equation}\label{intrans2}
\frac{L^{(1)}(({\bf x},1)({\bf y},1))}{L^{(m+1)}(({\bf x},1)({\bf y},1))}\mapsto
\frac{L^{(i)}((a_0^{(1)}x^{(1)}+x^{(2)},\ldots ,a_0^{(m-1)}x^{(1)}+x^{(m)},a_0^{(m)}x^{(1)}+1),({\bf y},1))}{L^{(m+1)}((a_0^{(1)}x^{(1)}+x^{(2)},\ldots ,a_0^{(m-1)}x^{(1)}+x^{(m)},a_0^{(m)}x^{(1)}+1),({\bf y},1))},
\end{equation}
for each $i\in\{1,\ldots m\}$. If we write 
$$L^{(i)}((x^{(1)},\ldots x^{(m+1)}),(y^{(1)},\ldots y^{(m+1)}))=\sum_{1\leq j,l\leq m+1} c_{j,l}^{(i)}x^{(j)}y^{(l)},\ i\in\{1,\ldots m+1\},$$
and $C^{(i)}=\left[c_{j,l}^{(i)}\right]_{1\leq j,l\leq m+1}$, then by the transformation \eqref{intrans2} we get the following transformation of the matrices $C^{(i)}$, $i\in\{1,\ldots m+1\}$:
$$C^{(i)}\mapsto\left[\begin{array}{ccccc}
a_0^{(1)} & a_0^{(2)} & \cdots & a_0^{(m)} & 1\\
1 & 0 & \cdots & 0 & 0\\
0 & 1 & \cdots & 0 & 0\\
\vdots & \vdots & \ddots & \vdots & \vdots\\
0 & 0 & \cdots & 1 & 0
\end{array}\right]C^{(i)},\ i\in\{1,\ldots m+1\}.$$

A similar transformation holds for $\mathbf y$ and instead of describing input transformation for the $m$-tuple ${\bf y}$ in detail it will be suffices to swap ${\bf x}$ with ${\bf y}$ in the algorithm. 

From the Jacobi-Perron algorithm \ref{JPA} we know that after input transformation for ${\bf x}$ and ${\bf y}$ these two $m$-tuples become $m$-tuples of formal symbols that attain only positive values. Therefore, for each $i\in\{1,\ldots ,m\}$ the value of $\frac{L^{(i)}(({\bf x},1)({\bf y},1))}{L^{(m+1)}(({\bf x},1)({\bf y},1))}$ lies between the minimum and maximum of the values $\frac{c_{j,l}^{(i)}}{c_{j,l}^{(m+1)}}$, $j,l\in\{1,\ldots ,m+1\}$, on condition that all the values $c_{j,l}^{(m+1)}$, $j,l\in\{1,\ldots ,m+1\}$, have the same sign. This fact, just like in the case of a single variable ${\bf x}$, comes from the mediant inequality \eqref{med} and induction with respect to the number of the quotients $\frac{c_{j,l}^{(i)}}{c_{j,l}^{(m+1)}}$. Thus, if all the mentioned values have the same integral part, then also $\frac{L^{(i)}(({\bf x},1)({\bf y},1))}{L^{(m+1)}(({\bf x},1)({\bf y},1))}$ has this value as the integral part. Then, we have $d^{(i)}=\left\lfloor\frac{L^{(i)}(({\bf x},1)({\bf y},1))}{L^{(m+1)}(({\bf x},1)({\bf y},1))}\right\rfloor=\left\lfloor\frac{c_{1,1}^{(i)}}{c_{1,1}^{(m+1)}}\right\rfloor$. Summing up, if $$\left\lfloor\frac{c_{j_1,l_1}^{(i)}}{c_{j_1,l_1}^{(m+1)}}\right\rfloor =\left\lfloor\frac{c_{j_2,l_2}^{(i)}}{c_{j_2,l_2}^{(m+1)}}\right\rfloor$$
for each $i\in\{1,\ldots m\}$ and $j_1,j_2,l_1,l_2\in\{1,\ldots m+1\}$, then we get the next terms of our output sequences ${\bf d}^{(i)}$, $i\in\{1,\ldots ,m\}$. Next, we transform
\begin{align*}
\frac{L^{(1)}(({\bf x},1),({\bf y},1))}{L^{(m+1)}(({\bf x},1),({\bf y},1))}&\mapsto\frac{1}{\frac{L^{(m)}(({\bf x},1),({\bf y},1))}{L^{(m+1)}(({\bf x},1),({\bf y},1))}-d^{(m)}}\\
&=\frac{L^{(m+1)}(({\bf x},1),({\bf y},1))}{L^{(m)}(({\bf x},1),({\bf y},1))-d^{(m)}L^{(m+1)}(({\bf x},1),({\bf y},1))},\\
\frac{L^{(i+1)}(({\bf x},1),({\bf y},1))}{L^{(m+1)}(({\bf x},1),({\bf y},1))}&\mapsto\frac{\frac{L^{(i)}(({\bf x},1),({\bf y},1))}{L^{(m+1)}(({\bf x},1),({\bf y},1))}-d^{(i)}}{\frac{L^{(m)}(x,y({\bf x},1),({\bf y},1))}{L^{(m+1)}(({\bf x},1),({\bf y},1))}-d^{(m)}}\\
&=\frac{L^{(i)}(({\bf x},1),({\bf y},1))-d^{(i)}L^{(m+1)}(({\bf x},1),({\bf y},1))}{L^{(m)}(({\bf x},1),({\bf y},1))-d^{(m)}L^{(m+1)}(({\bf x},1),({\bf y},1))},\, i\in\{1,\ldots m-1\}.
\end{align*}
Hence,
\begin{align*}
    C^{(1)}&\mapsto C^{(m+1)},\\
    C^{(i+1)}&\mapsto C^{(i)}-d^{(i)}C^{(m+1)},\, i\in\{1,\ldots m-1\},
\end{align*}
and
$$\textbf{d}^{(i)}=\left(\textbf{d}^{(i)},d^{(i)}\right), i\in\{1,\ldots m\}.$$

According to the above explanation, the Gosper's algorithm for ${\bf x}, {\bf y}\in\R^m$ and bilinear forms $C_1,\ldots ,C_{m+1}$ with integral coefficients as inputs can be described as in Algorithm \ref{Alg: GospMCF2}.

\newpage
\IncMargin{1.5em}
\begin{algorithm}[h]\label{Alg: GospMCF2}
	\caption{Algorithm for obtaining the MCF expansion of $\left(\frac{L^{(1)}(({\bf x},1)({\bf y},1))}{L^{(m+1)}(({\bf x},1)({\bf y},1))},\ldots ,\frac{L^{(m)}(({\bf x},1)({\bf y},1))}{L^{(m+1)}(({\bf x},1)({\bf y},1))}\right)$ given the MCF expansion of $\mathbf x = (x^{(1)}, \ldots, x^{(m)})$ and $\mathbf y = (y^{(1)}, \ldots, y^{(m)})$, where $L^{(1)}, \ldots, L^{(m+1)}$ are bilinear forms.}
	\SetKwData{Left}{left}
	\SetKwData{This}{this}
	\SetKw{And}{and}
	\SetKwFunction{Union}{Union}
	\SetKwFunction{FindCompress}{FindCompress}
	\SetKwInOut{Input}{Input}
	\SetKwInOut{Output}{Output}
    \SetKw{KwTo}{to}
    \SetKwComment{Comment}{$\triangleright$ }{ }
	\Input{$\mathbf x =[(a_n^{(1)})_{n\geq0},...,(a_n^{(m)})_{n\geq0}]$, $\mathbf y= [(b_n^{(1)})_{n\geq0},...,(b_n^{(m)})_{n\geq0}]$, $C^{(i)}=\left[c_{j,l}^{(i)}\right]_{1\leq j,l\leq m+1}\in M(m+1\times m+1,\Z)$, $i\in\{1,\ldots,m+1\}$, $M,N \in \mathbb N$}
	\Output{$[(d_n^{(1)})_{n\geq0},...,(d_n^{(m)})_{n\geq0}] = \left(\frac{L^{(1)}(({\bf x},1)({\bf y},1))}{L^{(m+1)}(({\bf x},1)({\bf y},1))},\ldots ,\frac{L^{(m)}(({\bf x},1)({\bf y},1))}{L^{(m+1)}(({\bf x},1)({\bf y},1))}\right)$}
	\BlankLine
	$r \gets 0$, $s \gets 0$, $t \gets 0$

    \While{$s<M$ \And $r<N$}{

        \If{$\sgn\left(c_{j_1,l_1}^{(m+1)}\right)=\sgn\left(c_{j_2,l_2}^{(m+1)}\right)$ and $\left\lfloor\frac{c_{j_1,l_1}^{(i)}}{c_{j_1,l_1}^{(m+1)}}\right\rfloor =\left\lfloor\frac{c_{j_2,l_2}^{(i)}}{c_{j_2,l_2}^{(m+1)}}\right\rfloor$ for each $j_1,j_2,l_1,l_2\in\{1,\ldots m+1\}$}{
            $d_s^{(i)} \gets\left\lfloor\frac{c_{1,1}^{(i)}}{c_{1,1}^{(m+1)}}\right\rfloor$ for each $i\in\{1,\ldots,m\}$
        
             $C^{(1)} \gets C^{(m+1)}$,\,
             $C^{(i+1)} \gets C^{(i)}-d^{(i)}_sC^{(m+1)},\, i\in\{1,\ldots m\}$

            $s \gets s + 1$
        }

        \Else{
        
        $C^{(i)}\gets\left[\begin{array}{ccccc}
        a_{\lfloor t/2\rfloor}^{(1)} & a_{\lfloor t/2\rfloor}^{(2)} & \cdots & a_{\lfloor t/2\rfloor}^{(m)} & 1\\
        1 & 0 & \cdots & 0 & 0\\
        0 & 1 & \cdots & 0 & 0\\
        \vdots & \vdots & \ddots & \vdots & \vdots\\
        0 & 0 & \cdots & 1 & 0
        \end{array}\right]C^{(i)},\ i\in\{1,\ldots m+1\},$

${\bf a}^{(i)}\leftrightarrow {\bf b}^{(i)},\, i\in\{1,\ldots m\}$

$C^{(i)}\gets \left(C^{(i)}\right)^T,\, i\in\{1,\ldots m+1\}$
        
$t \gets t+1$
            
        }
            
    $r \gets r+1$

    }
    
    \Return $[(d_n^{(1)})_{n=0}^{s-1},...,(d_n^{(m)})_{n=0}^{s-1}]$

\end{algorithm}
\DecMargin{1.5em}

Example \ref{E3} in Appendix illustrates the above algorithm.

\subsection{Proof of feasibility of the algorithm}\label{feas2}

We show that after a finite number of input transformations we can perform an output transformation.

We see that the $n_1$-fold iteration of the input transformation of the matrix $C^{(i)}$ with respect to ${\bf x}$ and $n_2$-fold iteration of the input transformation of $C$ with respect to ${\bf y}$ is
\begin{equation*}
\begin{split}
C^{(i)}\mapsto &\left(\prod_{k=0}^{n_1-1}\left[\begin{array}{cccccc}
a_k^{(1)} & 1 & \cdots & 0 & 0 \\
\vdots & \vdots & \ddots & \vdots & \vdots\\
a_k^{(m-1)} & 0 & \cdots & 1 & 0\\
a_k^{(m)} & 0 & \cdots & 0 & 1\\
1 & 0 & \cdots & 0 & 0
\end{array}\right]\right)^T\cdot C^{(i)}\cdot \prod_{k=0}^{n_2-1}\left[\begin{array}{cccccc}
b_k^{(1)} & 1 & \cdots & 0 & 0 \\
\vdots & \vdots & \ddots & \vdots & \vdots\\
b_k^{(m-1)} & 0 & \cdots & 1 & 0\\
b_k^{(m)} & 0 & \cdots & 0 & 1\\
1 & 0 & \cdots & 0 & 0
\end{array}\right]\\
&=[A_{n_1-j}^{(i)}]_{1\leq i,j\leq m+1}^T\cdot C\cdot [B_{n_2-l}^{(i)}]_{1\leq i,l\leq m+1},
\end{split}
\end{equation*}
where the values $B_l^{(i)}$ are defined by \eqref{convmat} for $\textbf{b}$ in the place of $\textbf{a}$. This means that meanwhile the bilinear forms $L^{(i)}$, $i\in\{1,\ldots ,m+1\}$, transform as follows:
\begin{equation*}
  \begin{split}
& L^{(i)}((x^{(1)},\ldots ,x^{(m+1)}),(y^{(1)},\ldots ,y^{(m+1)}))\mapsto \\
& L^{(i)}\left(\left(\sum_{j=1}^{m+1}A_{n_1+1-j}^{(1)}x^{(j)},\ldots ,\sum_{j=1}^{m+1}A_{n_1+1-j}^{(m+1)}x^{(j)}\right),\left(\sum_{j=1}^{m+1}B_{n_2+1-j}^{(1)}y^{(j)},\ldots ,\sum_{j=1}^{m+1}B_{n_2+1-j}^{(m+1)}y^{(j)}\right)\right)\\
&=\sum_{1\leq j,l\leq m+1}\left(\sum_{1\leq k,h\leq m+1} c_{k,h}^{(i)}A_{n_1+1-j}^{(k)}B_{n_2+1-l}^{(h)}\right)x^{(j)}y^{(l)},\, i\in\{1,\ldots m\}.
  \end{split}  
\end{equation*}
Hence,
\begin{equation}\label{coefftrans2}
    c_j^{(i)}\mapsto\sum_{1\leq k,h\leq m+1} c_{k,h}^{(i)}A_{n_1+1-j}^{(k)}B_{n_2+1-l}^{(h)},\, i,j,l\in\{1,\ldots ,m+1\}
\end{equation}
and
\begin{equation}\label{ratiotrans2}
  \begin{split}
&\frac{L^{(i)}(({\bf x},1),({\bf y},1))}{L^{(m+1)}(({\bf x},1),({\bf y},1))}\mapsto\\
&=\frac{\sum_{1\leq j,l\leq m+1}\left(\sum_{1\leq k,h\leq m+1} c_{k,h}^{(i)}A_{n_1+1-j}^{(k)}B_{n_2+1-l}^{(h)}\right)x^{(j)}y^{(l)}}{\sum_{1\leq j,l\leq m+1}\left(\sum_{1\leq k,h\leq m+1} c_{k,h}^{(m+1)}A_{n_1+1-j}^{(k)}B_{n_2+1-l}^{(h)}\right)x^{(j)}y^{(l)}},\, i\in\{1,\ldots m\},
  \end{split}  
\end{equation}
where $x^{(m+1)}=y^{(m+1)}=1$. As all the values of $x^{(j)}$ and $y^{(l)}$ are positive, the last expression in \eqref{ratiotrans2} lies between the minimum and maximum of the values 
\begin{equation}\label{coeffratiotrans2}
    \frac{\sum_{1\leq k,h\leq m+1} c_{k,h}^{(i)}A_{n_1+1-j}^{(k)}B_{n_2+1-l}^{(h)}}{\sum_{1\leq k,h\leq m+1} c_{k,h}^{(m+1)}A_{n_1+1-j}^{(k)}B_{n_2+1-l}^{(h)}}=\frac{\sum_{1\leq k,h\leq m+1} c_{k,h}^{(i)}\cdot\frac{A_{n_1+1-j}^{(k)}}{A_{n_1+1-j}^{(m+1)}}\cdot\frac{B_{n_2+1-l}^{(h)}}{B_{n_2+1-l}^{(m+1)}}}{\sum_{1\leq k,h\leq m+1} c_{k,h}^{(m+1)}\cdot\frac{A_{n_1+1-j}^{(k)}}{A_{n_1+1-j}^{(m+1)}}\cdot\frac{B_{n_2+1-l}^{(h)}}{B_{n_2+1-l}^{(m+1)}}},\, j,l\in\{1,\ldots ,m+1\},
\end{equation}
on condition that all of the values $\sum_{1\leq k,h\leq m+1} c_{k,h}^{(m+1)}A_{n_1+1-j}^{(k)}B_{n_2+1-l}^{(h)}$, $j,l\in\{1,\ldots ,m+1\}$, have the same sign. Since the Jacobi-Perron algorithm for ${\bf x}$ and ${\bf y}$ never stops, we know that $\lim_{n_1\to\infty}\frac{A_{n_1+1-j}^{(k)}}{A_{n_1+1-j}^{(m+1)}}=x_0^{(k)}$, $k\in\{1,\ldots m\}$, $j\in\{1,\ldots ,m+1\}$, and $\lim_{n_2\to\infty}\frac{B_{n_2+1-l}^{(h)}}{B_{n_2+1-l}^{(m+1)}}=y_0^{(l)}$, $h\in\{1,\ldots m\}$, $l\in\{1,\ldots ,m+1\}$. From this, \eqref{coefftrans2}, and \eqref{coeffratiotrans2} we see that the consecutive transformations of the ratios $\frac{c_{j,l}^{(i)}}{c_{j,l}^{(m+1)}}$, $j,l\in\{1,\ldots ,m+1\}$, converge to the same value, namely $$\frac{L^{(i)}((x_0^{(1)},\ldots ,x_0^{(m)},1),(y_0^{(1)},\ldots ,y_0^{(m)},1))}{L^{(m+1)}((x_0^{(1)},\ldots ,x_0^{(m)},1),(y_0^{(1)},\ldots ,y_0^{(m)},1))}.$$ 
In particular, for sufficiently large $n_1,n_2\in\N$ all the values 
$$\sum_{1\leq k,h\leq m+1} c_{k,h}^{(m+1)}\cdot\frac{A_{n_1+1-j}^{(k)}}{A_{n_1+1-j}^{(m+1)}}\cdot\frac{B_{n_2+1-l}^{(h)}}{B_{n_2+1-l}^{(m+1)}}, \quad j,l\in\{1,\ldots ,m+1\}$$ 
have the same sign. Consequently, all the values 
$$\sum_{1\leq k,h\leq m+1} c_{k,h}^{(m+1)}A_{n_1+1-j}^{(k)}B_{n_2+1-l}^{(h)}, \quad j,l\in\{1,\ldots ,m+1\}$$ 
have the same sign as all of the values $A_j^{(m+1}$, $B_l^{(m+1)}$, $j,l\in\N$, are positive. Since the $m+1$-tuple 
$$\left(\frac{L^{(1)}((x_0^{(1)},\ldots ,x_0^{(m)},1),(y_0^{(1)},\ldots ,y_0^{(m)},1))}{L^{(m+1)}((x_0^{(1)},\ldots ,x_0^{(m)},1),(y_0^{(1)},\ldots ,y_0^{(m)},1))},\ldots ,\frac{L^{(m)}((x_0^{(1)},\ldots ,x_0^{(m)},1),(y_0^{(1)},\ldots ,y_0^{(m)},1))}{L^{(m+1)}((x_0^{(1)},\ldots ,x_0^{(m)},1),(y_0^{(1)},\ldots ,y_0^{(m)},1))},1\right)$$ 
is $\Q$-linearly independent, for each $i\in\{1,\ldots ,m\}$ the value  
$$\frac{L^{(i)}((x_0^{(1)},\ldots ,x_0^{(m)},1),(y_0^{(1)},\ldots ,y_0^{(m)},1))}{L^{(m+1)}((x_0^{(1)},\ldots ,x_0^{(m)},1),(y_0^{(1)},\ldots ,y_0^{(m)},1))}$$ 
is irrational, in particular non-integral. This means that for sufficiently large $n_1$, $n_2$, after $n_1$ input transformations with respect to ${\bf x}$ and $n_2$ input transformations with respect to ${\bf y}$ the ratios $\frac{c_{j,l}^{(i)}}{c_{j,l}^{(m+1)}}$, $j,l\in\{1,\ldots ,m+1\}$ have the same integral part, so we may perform output transformation.

\section{Experimental results} \label{sec:exp}
In this section we collect some data about the performance of the algorithms. In particular, we examine the number of inputs that are required at each step before getting an output. The main outcome of our computations is that the number of inputs required by Algorithm \ref{Alg: GospMCF1} and \ref{Alg: GospMCF2} grows linearly with the number of outputs. Moreover, it seems to approach always the same line, independently of the transformation $C$ and of the sequence of partial quotients $((a_n)_{n\geq 0},(b_n)_{n\geq 0})$, both in the case of cubic irrationals and in the case of random sequences. The quite surprising outcome is that this line has slope less than $1$, around $\frac{1}{4}$ for Algorithm \ref{Alg: GospMCF1} and around $\frac{1}{2}$ for Algorithm \ref{Alg: GospMCF2}, independently of the dimension $m$ of the MCF. This means that we usually expect to get $4$ outputs by using only, respectively, $1$ input for Algorithm \ref{Alg: GospMCF1} and $2$ inputs for Algorithm \ref{Alg: GospMCF2}. The SageMath code used for this analysis is publicly available\footnote{\href{https://github.com/giulianoromeont/Continued-fractions}{https://github.com/giulianoromeont/Continued-fractions}}.

\begin{figure}[h]
\centering
\includegraphics[scale=0.6]{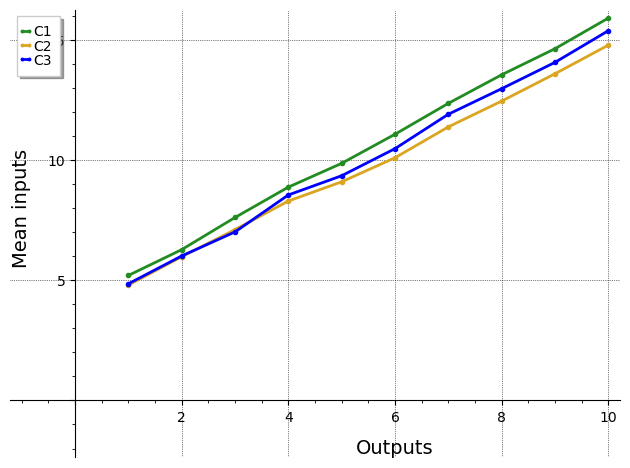}
\caption{ Mean number of inputs required by Algorithm \ref{Alg: GospMCF1} for $m=2$ to obtain outputs up to $10$, for $(\sqrt[3]{d},\sqrt[3]{d^2})$, $2\leq d\leq 100$ non-cube integer, using  $C_1$, $C_2$, $C_3$ as defined in \eqref{Eq: CCC}.}
\label{Fig: mean}
\end{figure}

In Figure \ref{Fig: mean}, we show the mean number of inputs required to get $10$ outputs with Algorithm \ref{Alg: GospMCF1}, $m=2$, for the transformations
\begin{equation}\label{Eq: CCC}
C_1=\left[\begin{array}{ccc}
            2 & 0 & 0 \\
            0 & 2 & 0 \\
            0 & 0 & 1
        \end{array}\right], \ \ \ \  \ \ C_2=\left[\begin{array}{ccc}
            1 & -1 & 0 \\
            1 & -1 & 0 \\
            0 & 0 & 1
        \end{array}\right], \ \ \ \  \ \ C_3=\left[\begin{array}{ccc}
            3 & 5 & 0 \\
            5 & 3 & 0 \\
            1 & 0 & 2
        \end{array}\right],
\end{equation}
applied to the pair $(\sqrt[3]{d},\sqrt[3]{d^2})$ for $2\leq d\leq 100$ non-cube integer, i.e. to compute the MCF of
\[(2\sqrt[3]{d},2\sqrt[3]{d^2}), \ \ \  (\sqrt[3]{d}+\sqrt[3]{d^2},\sqrt[3]{d}-\sqrt[3]{d^2}),\ \ \  \left(\frac{3\sqrt[3]{d}+5\sqrt[3]{d^2}}{\sqrt[3]{d}+2},\frac{5\sqrt[3]{d}+3\sqrt[3]{d^2}}{\sqrt[3]{d}+2}\right).\]

\begin{figure}[h]
\centering
\includegraphics[scale=0.55]{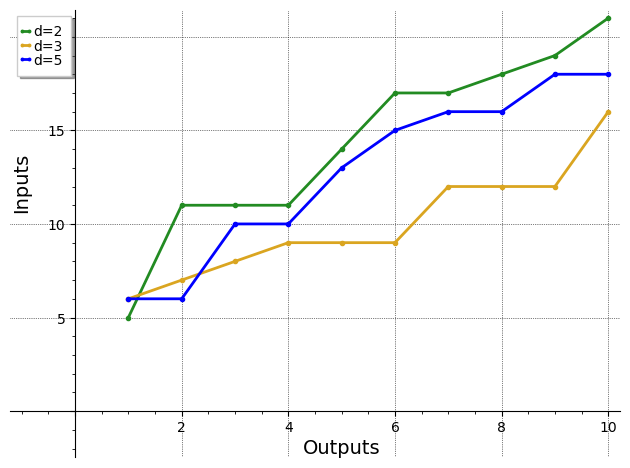}
\caption{Number of inputs required by Algorithm \ref{Alg: GospMCF1} for $m=2$ to obtain outputs up to $10$, for $(\sqrt[3]{d},\sqrt[3]{d^2})$, $d=2,3,5$, using $C_1$ as defined in \eqref{Eq: CCC}}
\label{Fig: d235} 
\end{figure}

It is possible to see that the mean number of inputs tends to have a similar trend for different transformations. The behavior of the total number of inputs required for each output is in general various, as shown in Figure \ref{Fig: d235} for $d=2,3,5$.

\begin{center}
\begin{figure}[!htb]
\minipage{0.15\textwidth}
  \includegraphics[scale=0.27]{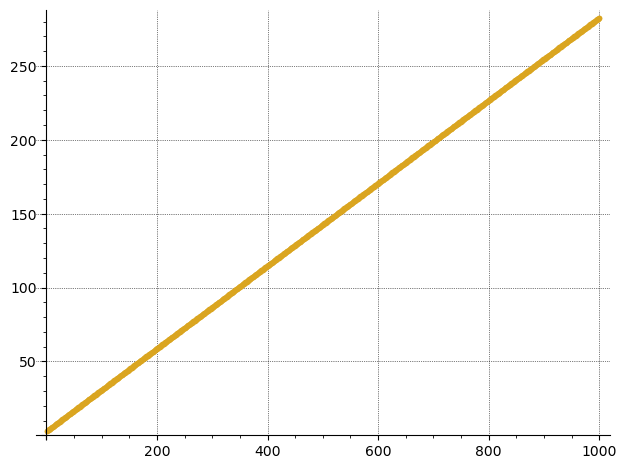}
\endminipage \hspace{2cm}
\minipage{0.15\textwidth}
  \includegraphics[scale=0.27]{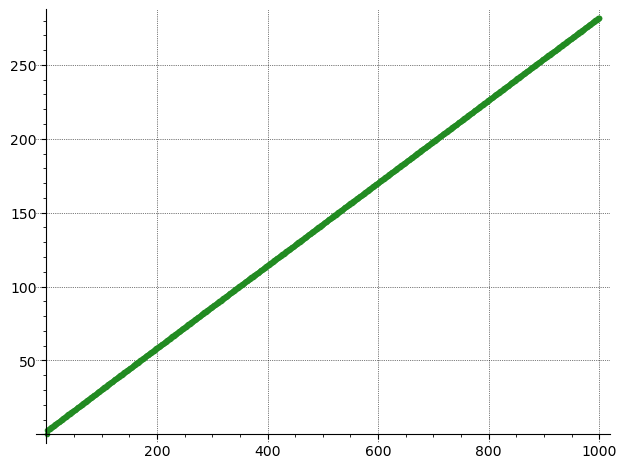}
\endminipage\vspace{0.5cm}\\
\minipage{0.15\textwidth}
  \includegraphics[scale=0.27]{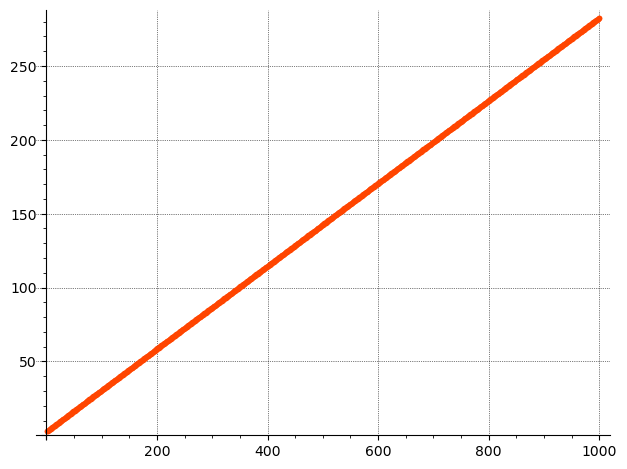}
\endminipage\hspace{2cm}
\minipage{0.15\textwidth}%
  \includegraphics[scale=0.27]{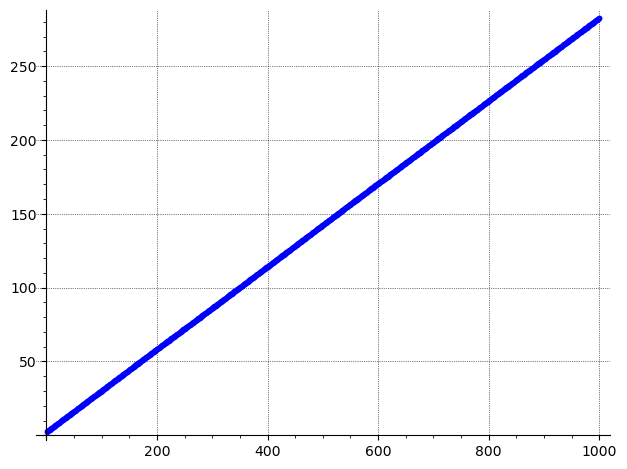}
\endminipage
\caption{Mean number of inputs required by Algorithm \ref{Alg: GospMCF1} for $m=2$ up to $1000$ outputs, for $4$ different $3\times 3$ matrices with entries in $\{0,\ldots,1000\}$. The mean is computed among $1000$ different random input sequences $((a_n),(b_n))$ with $b_n\leq a_n\leq 1000$.}
\label{Fig: 4graphics}
\end{figure}
\end{center}

In Figure \ref{Fig: 4graphics}, we plot the results for four different transformations $C$ applied to $1000$ different random sequences of partial quotients in $\{0,\ldots,10000\}$. The coefficients of $C$ are chosen randomly in $\{0,\ldots,10000\}$.

\begin{center}
\begin{figure}[!ht]
\includegraphics[scale=0.48]{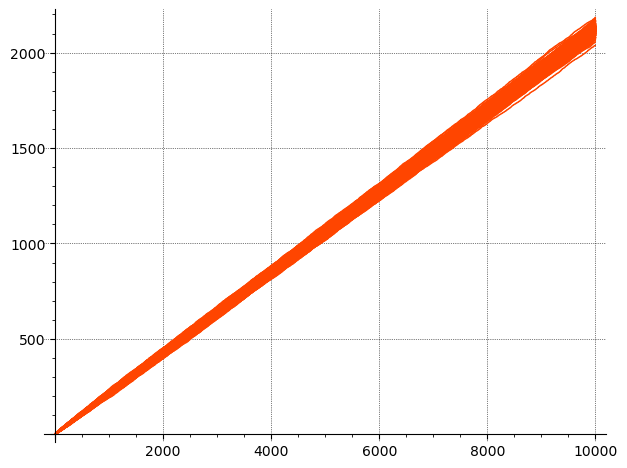}
\caption{Number of inputs required by Algorithm \ref{Alg: GospMCF1} for $m=2$ to get up to $10000$ outputs, for $1000$ different random sequences of partial quotients in $\{0,\ldots,10000\}$ and for $10000$ random transformations with coefficients in $\{0,\ldots,10000\}$.}
\label{Fig: Single10k}
\end{figure}
\end{center}

In Figure \ref{Fig: Single10k}, we compute the number of inputs for $1000$ different random sequences, and we plot them all together. Notice that they do not deviate so much from their mean.

\begin{center}
\begin{figure}[!htb]
\minipage{0.3\textwidth}
  \includegraphics[scale=0.42]{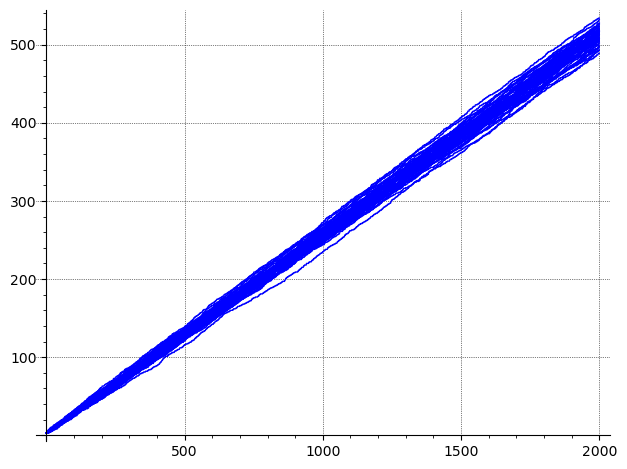}
\endminipage \hspace{2cm}
\minipage{0.3\textwidth}
  \includegraphics[scale=0.42]{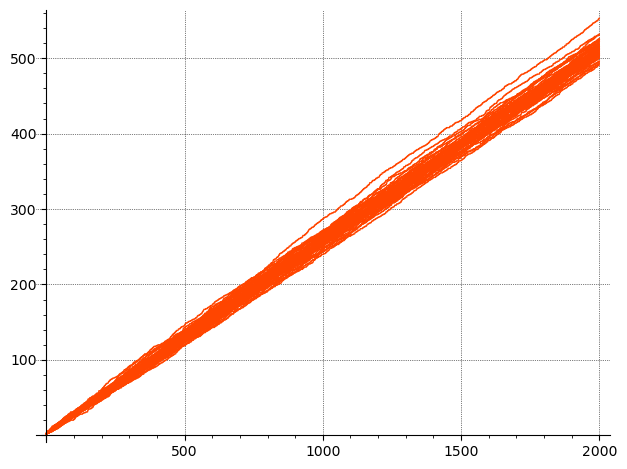}
\endminipage\vspace{0.5cm}
\caption{Number of inputs required by Algorithm \ref{Alg: GospMCF1} for $m=2$ to get up to $2000$ outputs for $500$ different random sequences. The partial quotients are randomly selected in $\{0,\ldots,10\}$ for the left graphic and in $\{0,\ldots,1000000\}$ for the right graphic.}
\label{Figure: diffsize}
\end{figure}
\end{center}

Moreover, it is possible to see in Figure \ref{Figure: diffsize} that the distribution around the same line does not depend on the size of the partial quotients. In fact, the number of inputs required to provide up to $2000$ outputs it is very similar in both the cases where the partial quotients are selected at random inside $\{0,\ldots,10\}$ or inside $\{0,\ldots,1000000\}$. The same result can be observed if we vary the size of the coefficients of the transformation matrix $C$.

Notice that, in all these cases, the number of inputs distributes around a line that has slope less than $1$, around $\frac{1}{4}$. This means that to get $N$ partial quotients of a transformation of $(\alpha,\beta)$, we expect to use around $\frac{N}{4}$ partial quotients of their continued fractions.\bigskip

\begin{center}
\begin{figure}[!ht]
\includegraphics[scale=0.45]{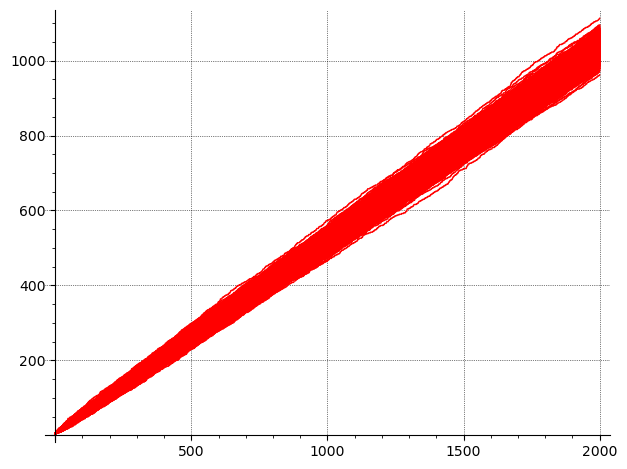}
\caption{ Number of inputs required by Algorithm \ref{Alg: GospMCF2} for $m=2$ to get up to $2000$ outputs, for the bilinear transformations using $2000$ different pairs of random sequences of partial quotients in $\{0,\ldots,2000\}$ and for $2000$ random transformations $C_1,C_2,C_3$, with coefficients in $\{0,\ldots,2000\}$.}
\label{Fig: biquad}
\end{figure}
\end{center}

In Figure \ref{Fig: biquad}, we show that a similar behavior is observed also in the case of bilinear transformations, using Algorithm \ref{Alg: GospMCF2}. Notice that the difference between Figure \ref{Fig: Single10k}, for Algorithm \ref{Alg: GospMCF1}, and Figure \ref{Fig: biquad}, for Algorithm \ref{Alg: GospMCF2}, is the slope of the line. The slope of the line for all sequences using Algorithm \ref{Alg: GospMCF2} is less than $1$ and it appears to be slightly larger than $\frac{1}{2}$.

Now we compare these results with the case $m=1$ and $m=3$, i.e. for unidimensional continued fractions and for the Jacobi-Perron algorithm working on triples of real numbers. It turns out that the results are similar and the number of inputs required grows linearly with the number of outputs obtained.

\begin{center}
\begin{figure}[!htb]
\minipage{0.3\textwidth}
  \includegraphics[scale=0.42]{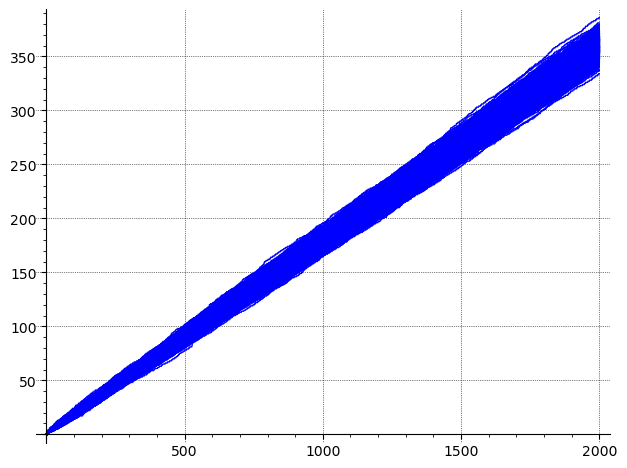}
\endminipage \hspace{2cm}
\minipage{0.3\textwidth}
  \includegraphics[scale=0.42]{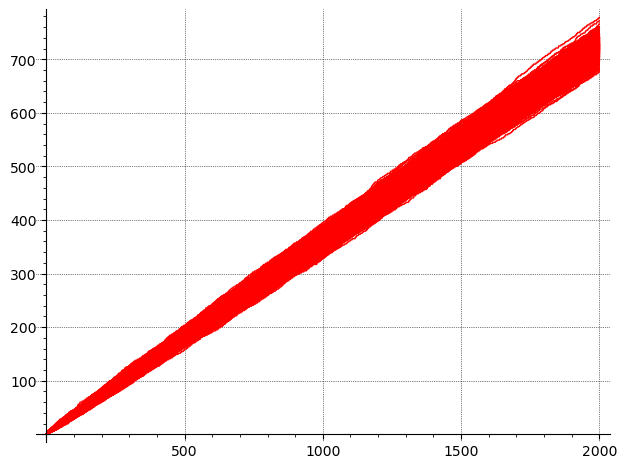}
\endminipage\vspace{0.5cm}
\caption{Number of inputs required for unidimensional continued fractions, $m=1$, to get up to $2000$ outputs for $2000$ different random sequences, with partial quotients randomly selected in $\{0,\ldots,2000\}$, for the single transformation of Algorithm \ref{Alg: Gosp1} (on the left) and the bilinear transformation of Algorithm \ref{Alg: Gosp2} (on the right).}
\label{Figure: unidim}
\end{figure}
\end{center}

\begin{center}
\begin{figure}[!htb]
\minipage{0.3\textwidth}
  \includegraphics[scale=0.42]{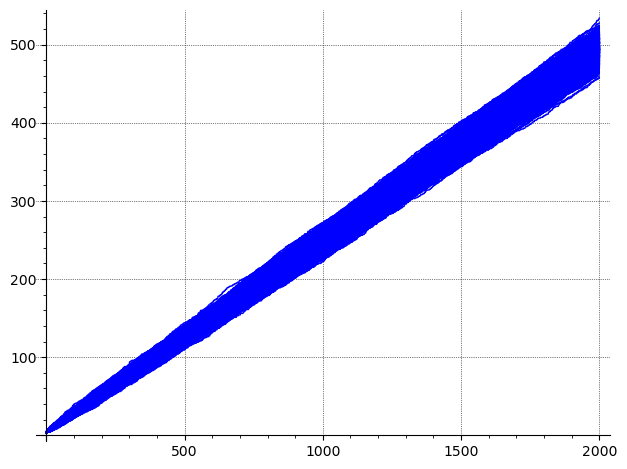}
\endminipage \hspace{2cm}
\minipage{0.3\textwidth}
  \includegraphics[scale=0.42]{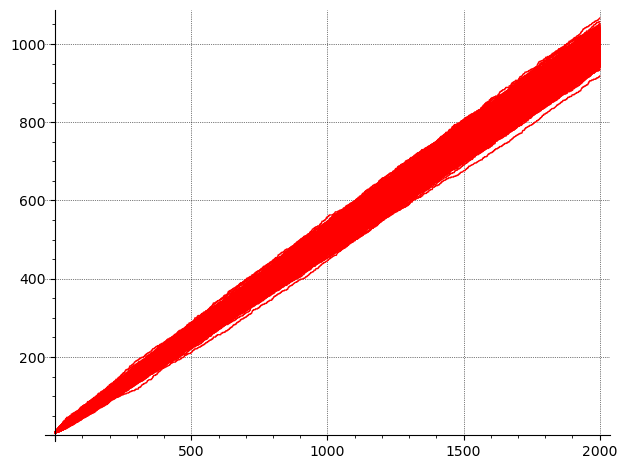}
\endminipage\vspace{0.5cm}
\caption{Number of inputs required for $3$-dimensional continued fractions, $m=3$, to get up to $2000$ outputs for $2000$ different random sequences, with partial quotients randomly selected in $\{0,\ldots,2000\}$, for the single transformation of Algorithm \ref{Alg: GospMCF1} (on the left) and the bilinear transformation of Algorithm \ref{Alg: GospMCF2} (on the right).}
\label{Figure: tridim}
\end{figure}
\end{center}

In Figure \ref{Figure: unidim} and Figure \ref{Figure: tridim}, we plot the results for both Algorithm \ref{Alg: GospMCF1} and the bilinear transformation of Algorithm \ref{Alg: GospMCF2} for, respectively, unidimensional continued fractions, i.e. $m=1$, and MCF with $m=3$. Notice that the results for Algorithm \ref{Alg: GospMCF1} for $m=2$, in Figure \ref{Figure: diffsize}, and $m=3$, in Figure \ref{Figure: tridim}, are really comparable, i.e. the number of inputs increase linearly and arrives to around $500$ inputs to get $2000$ inputs. For the unidimensional case of Algorithm \ref{Alg: Gosp1} the number of inputs increase linearly and arrives to around $350$ inputs to get $2000$ outputs. Also for the bilinear transformation of Algorithm \ref{Alg: GospMCF2}, the results for $m=2$ and $m=3$ are similar and around $1000$ inputs are required for $2000$ outputs. In the classical unidimensional case of Algorithm \ref{Alg: Gosp2}, around $700$ inputs are required for $2000$ outputs.

\section{Possible modification of the algorithm - partial output} \label{sec:mod}

After the first input transformations we know that our input $m$-tuple ${\bf x}$ is a string of variables that take on only positive values. This means that $\frac{L^{(i)}({\bf x},1)}{L^{(m+1)}({\bf x},1)}$ lies between the values of $\frac{c_j^{(i)}}{c_j^{(m+1)}}$, $j\in\{1,\ldots ,m+1\}$, on condition that all the values of $c_j^{(m+1)}$, $j\in\{1,\ldots ,m+1\}$, have the same sign (in this case we may assume that $c_j^{(m+1)}>0$, $j\in\{1,\ldots ,m+1\}$; we change the sign of all the values of $c_j^{(i)}$, if necessary). Hence, we may subtract 
\begin{align*}
    \tilde{b}^{(i)}=\min_{1\leq j\leq m+1}\left\lfloor\frac{c_j^{(i)}}{c_j^{(m+1)}}\right\rfloor ,\, i\in\{1,\ldots ,m\}
\end{align*} 
from $\frac{L^{(i)}({\bf x},1)}{L^{(m+1)}({\bf x},1)}$ for each $i\in\{1,\ldots ,m\}$. Then $\frac{L^{(i)}({\bf x},1)}{L^{(m+1)}({\bf x},1)}$ transforms in the following way.
\begin{align*}
    \frac{L^{(i)}({\bf x},1)}{L^{(m+1)}({\bf x},1)}\mapsto \frac{L^{(i)}({\bf x},1)}{L^{(m+1)}({\bf x},1)}-\tilde{b}^{(i)}=\frac{\sum_{j=1}^m (c_j^{(i)}-\tilde{b}^{(i)}c_j^{(m+1)})x^{(i)}+(c_{m+1}^{(i)}-\tilde{b}^{(i)}c_{m+1}^{(m+1)})}{\sum_{j=1}^m c_j^{(m+1)}x^{(i)}+c_{m+1}^{(m+1)}}
\end{align*}
This corresponds to the matrix transformation
\begin{align}\label{partoutmat}
    C\mapsto\left[\begin{array}{cccc}
        1 & \cdots & 0 & -\tilde{b}^{(1)} \\
        \vdots & \ddots & \vdots & \vdots \\
        0 & \cdots & 1 & -\tilde{b}^{(m)} \\
        0 & \cdots & 0 & 1
    \end{array}\right]\cdot C.
\end{align}
As the values $\tilde{b}^{(i)}$, $i\in\{1,\ldots ,m\}$, are only parts (with respect to addition) of our outputs ${b}^{(i)}$, $i\in\{1,\ldots ,m\}$, we call them the partial inputs and the operations
$$b^{(i)}\mapsto b^{(i)}+\tilde{b}^{(i)},\, i\in\{1,\ldots ,m\}$$
together with \eqref{partoutmat} are called the partial output transformation. This kind of transformation allows us not only to obtain parts of the output after smaller number of input transformations but also reduces the size of the values of $c_j^{(i)}$, $i\in\{1,\ldots ,m\}$, $j\in\{1,\ldots ,m+1\}$.

Example \ref{E2} in Appendix, compared with an appropriate part of Example \ref{E1}, shows the benefit of partial output.

We may do the same in case of two input $m$-tuples. Then, if only $\textbf{x}$, $\textbf{y}$ are the strings of positive numbers (which is true after the first input transformation with respect to both $\textbf{x}$ and $\textbf{y}$) and $C^{(m+1)}$ has all the entries positive, we may do the following partial output:
\begin{align*}
\frac{L^{(i)}(({\bf x},1),({\bf y},1))}{L^{(m+1)}(({\bf x},1),({\bf y},1))}&\mapsto\frac{L^{(i)}(({\bf x},1),({\bf y},1))}{L^{(m+1)}(({\bf x},1),({\bf y},1))}-\tilde{d}^{(i)}\\
&=\frac{L^{(i)}(({\bf x},1),({\bf y},1))-\tilde{d}^{(i)}L^{(m+1)}(({\bf x},1),({\bf y},1))}{L^{(m+1)}(({\bf x},1),({\bf y},1))},\, i\in\{1,\ldots m\},
\end{align*}
where $\tilde{d}^{(i)}=\min_{1\leq j,l\leq m+1} \left\lfloor\frac{c_{j,l}^{(i)}}{c_{j,l}^{(m+1)}}\right\rfloor$, $i\in\{1,\ldots m\}$. This can be implemented as
\begin{align*}
    \tilde{d}^{(i)}&\mapsto\min_{1\leq j,l\leq m_+1} \left\lfloor\frac{c_{j,l}^{(i)}}{c_{j,l}^{(m+1)}}\right\rfloor ,\,i\in\{1,\ldots m\}\\
    C^{(i)}&\mapsto C^{(i)}-\tilde{d}^{(i)}C^{(m+1)},\, i\in\{1,\ldots m\},\\
    {d}^{(i)}&\mapsto {d}^{(i)}+\tilde{d}^{(i)},\, i\in\{1,\ldots m\},
\end{align*}
where after each usual output transformation we need to put $d^{(i)}\mapsto 0$, $i\in\{1,\ldots m\}$. 

\section*{Acknowledgments}
The research of the first author was supported by  a grant of the National Science Centre (NCN), Poland, no. UMO-2019/34/E/ST1/00094. \\
The second and the third author are members of GNSAGA of INdAM.\\
G.R. acknowledges that this study was carried out within the MICS (Made in Italy – Circular and Sustainable) Extended Partnership and received funding from the European Union Next-Generation EU (PIANO NAZIONALE DI RIPRESA E RESILIENZA (PNRR) – MISSIONE 4 COMPONENTE 2, INVESTIMENTO 1.3 – D.D. 1551.11-10-2022, PE00000004).

\section*{Appendix}
In this Appendix, we see the behavior of the algorithms for MCFs over three examples.

\begin{Example}\label{E1}
    Let $(x^{(1)},x^{(2)})$ have the continued fraction expansion $[[1,\overline{1,2}],[\overline{1,0}]]$. Let us compute the beginning terms of the continued fraction expansion of $\left(\frac{x^{(1)}}{2},-\frac{x^{(2)}}{3}\right)$. Hence, our input is
    \begin{itemize}
        \item ${\bf a}^{(1)}=(1,\overline{1,2})$,
        \item ${\bf a}^{(2)}=(\overline{1,0})$,
        \item $C=\left[\begin{array}{ccc}
            3 & 0 & 0 \\
            0 & -2 & 0 \\
            0 & 0 & 6
        \end{array}\right]$,
    \end{itemize}
    whereas our output is $({\bf b}^{(1)},{\bf b}^{(2)})$.
    
    According to the algorithm, we start with the input transformation.
    \begin{align*}
        C\leftarrow &\left[\begin{array}{ccc}
            3 & 0 & 0 \\
            0 & -2 & 0 \\
            0 & 0 & 6
        \end{array}\right]\cdot\left[\begin{array}{ccc}
            1 & 1 & 0 \\
            1 & 0 & 1 \\
            1 & 0 & 0
        \end{array}\right]=\left[\begin{array}{ccc}
            3 & 3 & 0 \\
            -2 & 0 & -2 \\
            6 & 0 & 0
        \end{array}\right]\\
        {\bf a}^{(1)}\leftarrow &(\overline{1,2}),\,
        {\bf a}^{(2)}\leftarrow (\overline{0,1})
    \end{align*}
    The conditions $\left\lfloor\frac{c_1^{(1)}}{c_1^{(3)}}\right\rfloor =\left\lfloor\frac{c_2^{(1)}}{c_2^{(3)}}\right\rfloor =\left\lfloor\frac{c_3^{(1)}}{c_3^{(3)}}\right\rfloor$ and $\left\lfloor\frac{c_1^{(2)}}{c_1^{(3)}}\right\rfloor =\left\lfloor\frac{c_2^{(2)}}{c_2^{(3)}}\right\rfloor =\left\lfloor\frac{c_3^{(2)}}{c_3^{(3)}}\right\rfloor$ are not satisfied (the quotients $\frac{c_2^{(1)}}{c_2^{(3)}}$, $\frac{c_3^{(1)}}{c_3^{(3)}}$, $\frac{c_2^{(2)}}{c_2^{(3)}}$, $\frac{c_3^{(2)}}{c_3^{(3)}}$ are not even defined), so we perform the input transformation.
    \begin{align*}
        C\leftarrow &\left[\begin{array}{ccc}
            3 & 3 & 0 \\
            -2 & 0 & -2 \\
            6 & 0 & 0
        \end{array}\right]\cdot\left[\begin{array}{ccc}
            1 & 1 & 0 \\
            0 & 0 & 1 \\
            1 & 0 & 0
        \end{array}\right]=\left[\begin{array}{ccc}
            3 & 3 & 3 \\
            -4 & -2 & 0 \\
            6 & 6 & 0
        \end{array}\right]\\
        {\bf a}^{(1)}\leftarrow &(\overline{2,1}),\,
        {\bf a}^{(2)}\leftarrow (\overline{1,0})
    \end{align*}
    Still the conditions $\left\lfloor\frac{c_1^{(1)}}{c_1^{(3)}}\right\rfloor =\left\lfloor\frac{c_2^{(1)}}{c_2^{(3)}}\right\rfloor =\left\lfloor\frac{c_3^{(1)}}{c_3^{(3)}}\right\rfloor$ and $\left\lfloor\frac{c_1^{(2)}}{c_1^{(3)}}\right\rfloor =\left\lfloor\frac{c_2^{(2)}}{c_2^{(3)}}\right\rfloor =\left\lfloor\frac{c_3^{(2)}}{c_3^{(3)}}\right\rfloor$ are not satisfied (the quotients $\frac{c_3^{(1)}}{c_3^{(3)}}$, $\frac{c_3^{(2)}}{c_3^{(3)}}$ still are not defined), so we perform the input transformation once again.
    \begin{align*}
        C\leftarrow &\left[\begin{array}{ccc}
            3 & 3 & 3 \\
            -4 & -2 & 0 \\
            6 & 6 & 0
        \end{array}\right]\cdot\left[\begin{array}{ccc}
            2 & 1 & 0 \\
            1 & 0 & 1 \\
            1 & 0 & 0
        \end{array}\right]=\left[\begin{array}{ccc}
            12 & 3 & 3 \\
            -10 & -4 & -2 \\
            18 & 6 & 6 
        \end{array}\right]\\
        {\bf a}^{(1)}\leftarrow &(\overline{1,2}),\,
        {\bf a}^{(2)}\leftarrow (\overline{0,1})
    \end{align*}
    Now $\left\lfloor\frac{c_1^{(1)}}{c_1^{(3)}}\right\rfloor =\left\lfloor\frac{c_2^{(1)}}{c_2^{(3)}}\right\rfloor =\left\lfloor\frac{c_3^{(1)}}{c_3^{(3)}}\right\rfloor =0$ and $\left\lfloor\frac{c_1^{(2)}}{c_1^{(3)}}\right\rfloor =\left\lfloor\frac{c_2^{(2)}}{c_2^{(3)}}\right\rfloor =\left\lfloor\frac{c_3^{(2)}}{c_3^{(3)}}\right\rfloor =-1$, so we execute the output transformation.
    \begin{align*}
        C\leftarrow &\left[\begin{array}{ccc}
            0 & 0 & 1 \\
            1 & 0 & 0 \\
            0 & 1 & 1
        \end{array}\right]\cdot\left[\begin{array}{ccc}
            12 & 3 & 3 \\
            -10 & -4 & -2 \\
            18 & 6 & 6 
        \end{array}\right]=\left[\begin{array}{ccc}
            18 & 6 & 6 \\
            12 & 3 & 3 \\
            8 & 2 & 4 \\
            \end{array}\right]\\
        {\bf b}^{(1)}\leftarrow &(0),\,
        {\bf b}^{(2)}\leftarrow (-1)
    \end{align*}
    The conditions $\left\lfloor\frac{c_1^{(1)}}{c_1^{(3)}}\right\rfloor =\left\lfloor\frac{c_2^{(1)}}{c_2^{(3)}}\right\rfloor =\left\lfloor\frac{c_3^{(1)}}{c_3^{(3)}}\right\rfloor$ and $\left\lfloor\frac{c_1^{(2)}}{c_1^{(3)}}\right\rfloor =\left\lfloor\frac{c_2^{(2)}}{c_2^{(3)}}\right\rfloor =\left\lfloor\frac{c_3^{(2)}}{c_3^{(3)}}\right\rfloor$ are not satisfied, so we perform the input transformation.
    \begin{align*}
        C\leftarrow &\left[\begin{array}{ccc}
            18 & 6 & 6 \\
            12 & 3 & 3 \\
            8 & 2 & 4 \\
            \end{array}\right]\cdot\left[\begin{array}{ccc}
            1 & 1 & 0 \\
            0 & 0 & 1 \\
            1 & 0 & 0
        \end{array}\right]=\left[\begin{array}{ccc}
            24 & 18 & 6 \\
            15 & 12 & 3 \\
            12 & 8 & 2 \\
            \end{array}\right]\\
        {\bf a}^{(1)}\leftarrow &(\overline{2,1}),\,
        {\bf a}^{(2)}\leftarrow (\overline{1,0})
    \end{align*}
    The conditions $\left\lfloor\frac{c_1^{(1)}}{c_1^{(3)}}\right\rfloor =\left\lfloor\frac{c_2^{(1)}}{c_2^{(3)}}\right\rfloor =\left\lfloor\frac{c_3^{(1)}}{c_3^{(3)}}\right\rfloor$ and $\left\lfloor\frac{c_1^{(2)}}{c_1^{(3)}}\right\rfloor =\left\lfloor\frac{c_2^{(2)}}{c_2^{(3)}}\right\rfloor =\left\lfloor\frac{c_3^{(2)}}{c_3^{(3)}}\right\rfloor$ still are not satisfied, so we perform the input transformation once again.
    \begin{align*}
        C\leftarrow &\left[\begin{array}{ccc}
            24 & 18 & 6 \\
            15 & 12 & 3 \\
            12 & 8 & 2 \\
            \end{array}\right]\cdot\left[\begin{array}{ccc}
            2 & 1 & 0 \\
            1 & 0 & 1 \\
            1 & 0 & 0
        \end{array}\right]=\left[\begin{array}{ccc}
            72 & 24 & 18 \\
            45 & 15 & 12 \\
            34 & 12 & 8 \\
            \end{array}\right]\\
        {\bf a}^{(1)}\leftarrow &(\overline{1,2}),\,
        {\bf a}^{(2)}\leftarrow (\overline{0,1})
    \end{align*}
    Now $\left\lfloor\frac{c_1^{(1)}}{c_1^{(3)}}\right\rfloor =\left\lfloor\frac{c_2^{(1)}}{c_2^{(3)}}\right\rfloor =\left\lfloor\frac{c_3^{(1)}}{c_3^{(3)}}\right\rfloor =2$ and $\left\lfloor\frac{c_1^{(2)}}{c_1^{(3)}}\right\rfloor =\left\lfloor\frac{c_2^{(2)}}{c_2^{(3)}}\right\rfloor =\left\lfloor\frac{c_3^{(2)}}{c_3^{(3)}}\right\rfloor =1$, so we execute the output transformation.
    \begin{align*}
        C\leftarrow &\left[\begin{array}{ccc}
            0 & 0 & 1 \\
            1 & 0 & -2 \\
            0 & 1 & -1
        \end{array}\right]\cdot\left[\begin{array}{ccc}
            72 & 24 & 18 \\
            45 & 15 & 12 \\
            34 & 12 & 8 \\
            \end{array}\right]=\left[\begin{array}{ccc}
            34 & 12 & 8 \\
            4 & 0 & 2 \\
            11 & 3 & 4 \\
            \end{array}\right]\\
        {\bf b}^{(1)}\leftarrow &(0,2),\,
        {\bf b}^{(2)}\leftarrow (-1,1)
    \end{align*}

    Finally, by keeping performing the algorithm, we get that the continued fraction expansion for $\left(\frac{x^{(1)}}{2},-\frac{x^{(2)}}{3}\right)$ is $$[[0,2,2,2,1,5,2,\ldots ],[-1,1,0,2,1,3,1,\ldots ]].$$
\end{Example}

\begin{Example}\label{E2}
    Let us study once again the part of Example \ref{E1} after the first output transformation, when we try to get the second output. Now we are permitted to execute the partial output transformation. In this case we need to set $b^{(i)}\mapsto 0$, $i\in\{1,\ldots ,m\}$, after every (usual) output transformation. The set of our input data after first output is the following:
\begin{itemize}
    \item $C=\left[\begin{array}{ccc}
            18 & 6 & 6 \\
            12 & 3 & 3 \\
            8 & 2 & 4 \\
            \end{array}\right]$,
    \item ${\bf a}^{(1)}=(\overline{1,2})$,
    \item ${\bf a}^{(2)}=(\overline{0,1})$,
    \item ${\bf b}^{(1)}=(0)$,
    \item ${\bf b}^{(2)}=(-1)$.
\end{itemize}
We set $b^{(1)}=b^{(2)}=0$. Since $\left\lfloor\frac{18}{8}\right\rfloor =\left\lfloor\frac{6}{2}\right\rfloor\neq\left\lfloor\frac{6}{4}\right\rfloor$ and $\left\lfloor\frac{12}{8}\right\rfloor\neq\left\lfloor\frac{3}{2}\right\rfloor\neq\left\lfloor\frac{3}{4}\right\rfloor$, we compute
\begin{align*}
    \tilde{b}^{(1)}\leftarrow &\min\left\{\left\lfloor\frac{18}{8}\right\rfloor,\left\lfloor\frac{6}{2}\right\rfloor,\left\lfloor\frac{6}{4}\right\rfloor\right\}=1,\\
    \tilde{b}^{(2)}\leftarrow &\min\left\{\left\lfloor\frac{12}{8}\right\rfloor,\left\lfloor\frac{3}{2}\right\rfloor,\left\lfloor\frac{3}{4}\right\rfloor\right\}=0.
\end{align*}
Since at least one of the values of $\tilde{b}^{(1)}$, $\tilde{b}^{(2)}$ is nonzero, it makes sense to execute partial output transformation.
\begin{align*}
    C\leftarrow &\left[\begin{array}{ccc}
            1 & 0 & -1 \\
            0 & 1 & 0 \\
            0 & 0 & 1 \\
            \end{array}\right]\cdot\left[\begin{array}{ccc}
            18 & 6 & 6 \\
            12 & 3 & 3 \\
            8 & 2 & 4 \\
            \end{array}\right]=\left[\begin{array}{ccc}
            10 & 4 & 2 \\
            12 & 3 & 3 \\
            8 & 2 & 4 \\
            \end{array}\right]\\
    b^{(1)}\leftarrow &b^{(1)}+\tilde{b}^{(1)}=0+1=1\\
    b^{(2)}\leftarrow &b^{(2)}+\tilde{b}^{(2)}=0+0=0
\end{align*}
Now we do input transformation.
\begin{align*}
    C\leftarrow &\left[\begin{array}{ccc}
            10 & 4 & 2 \\
            12 & 3 & 3 \\
            8 & 2 & 4 \\
            \end{array}\right]\cdot\left[\begin{array}{ccc}
            1 & 1 & 0 \\
            0 & 0 & 1 \\
            1 & 0 & 0 \\
            \end{array}\right]=\left[\begin{array}{ccc}
            12 & 10 & 4 \\
            15 & 12 & 3 \\
            12 & 8 & 2 \\
            \end{array}\right]\\
    a^{(1)}\leftarrow &(\overline{2,1})\\
    a^{(2)}\leftarrow &(\overline{1,0})
\end{align*}
We have $\left\lfloor\frac{12}{12}\right\rfloor =\left\lfloor\frac{10}{8}\right\rfloor 
 \neq\left\lfloor\frac{4}{2}\right\rfloor$ and $\left\lfloor\frac{15}{12}\right\rfloor =\left\lfloor\frac{12}{8}\right\rfloor =\left\lfloor\frac{3}{2}\right\rfloor$. Compute
\begin{align*}
    \tilde{b}^{(1)}\leftarrow &\min\left\{\left\lfloor\frac{12}{12}\right\rfloor ,\left\lfloor\frac{10}{8}\right\rfloor ,\left\lfloor\frac{4}{2}\right\rfloor\right\}=1,\\
    \tilde{b}^{(2)}\leftarrow &\min\left\{\left\lfloor\frac{15}{12}\right\rfloor ,\left\lfloor\frac{12}{8}\right\rfloor ,\left\lfloor\frac{3}{2}\right\rfloor\right\}=1,
\end{align*}
so we perform partial output transformation.
\begin{align*}
    C\leftarrow &\left[\begin{array}{ccc}
            1 & 0 & -1 \\
            0 & 1 & -1 \\
            0 & 0 & 1 \\
            \end{array}\right]\cdot\left[\begin{array}{ccc}
            12 & 10 & 4 \\
            15 & 12 & 3 \\
            12 & 8 & 2 \\
            \end{array}\right]=\left[\begin{array}{ccc}
            0 & 2 & 2 \\
            3 & 4 & 1 \\
            12 & 8 & 2 \\
            \end{array}\right]\\
    b^{(1)}\leftarrow &b^{(1)}+\tilde{b}^{(1)}=1+1=2\\
    b^{(2)}\leftarrow &b^{(2)}+\tilde{b}^{(2)}=0+1=1
\end{align*}
Perform input transformation.
\begin{align*}
    C\leftarrow &\left[\begin{array}{ccc}
            0 & 2 & 2 \\
            3 & 4 & 1 \\
            12 & 8 & 2 \\
            \end{array}\right]\cdot\left[\begin{array}{ccc}
            2 & 1 & 0 \\
            1 & 0 & 1 \\
            1 & 0 & 0 \\
            \end{array}\right]=\left[\begin{array}{ccc}
            4 & 0 & 2 \\
            11 & 3 & 4 \\
            34 & 12 & 8 \\
            \end{array}\right]\\
    a^{(1)}\leftarrow &(\overline{1,2})\\
    a^{(2)}\leftarrow &(\overline{0,1})
\end{align*}
We have $\left\lfloor\frac{4}{34}\right\rfloor =\left\lfloor\frac{0}{12}\right\rfloor =\left\lfloor\frac{2}{8}\right\rfloor =0$ and $\left\lfloor\frac{11}{34}\right\rfloor =\left\lfloor\frac{3}{12}\right\rfloor =\left\lfloor\frac{4}{8}\right\rfloor =0$, so finally perform output transformation.
\begin{align*}
    \tilde{b}^{(1)}\leftarrow &\min\left\{\left\lfloor\frac{4}{34}\right\rfloor ,\left\lfloor\frac{0}{12}\right\rfloor ,\left\lfloor\frac{2}{8}\right\rfloor\right\}=0\\
    \tilde{b}^{(2)}\leftarrow &\min\left\{\left\lfloor\frac{11}{34}\right\rfloor ,\left\lfloor\frac{3}{12}\right\rfloor ,\left\lfloor\frac{4}{8}\right\rfloor\right\}=0\\
    C\leftarrow &\left[\begin{array}{ccc}
            0 & 0 & 1 \\
            1 & 0 & 0 \\
            0 & 1 & 0 \\
            \end{array}\right]\cdot\left[\begin{array}{ccc}
            4 & 0 & 2 \\
            11 & 3 & 4 \\
            34 & 12 & 8 \\
            \end{array}\right]=\left[\begin{array}{ccc}
            34 & 12 & 8 \\
            4 & 0 & 2 \\
            11 & 3 & 4 \\
            \end{array}\right]\\
    b^{(1)}\leftarrow &b^{(1)}+\tilde{b}^{(1)}=2+0\\
    b^{(2)}\leftarrow &b^{(2)}+\tilde{b}^{(2)}=1+0\\
    {\bf b}^{(1)}\leftarrow &({\bf b}^{(1)},b^{(1)})=(0,2)\\
    {\bf b}^{(2)}\leftarrow &({\bf b}^{(2)},b^{(2)})=(-1,1)
\end{align*}
Finally, our set of data after second output is just like in Example \ref{E1} without partial output transformation permitted:
\begin{itemize}
    \item $C=\left[\begin{array}{ccc}
            34 & 12 & 8 \\
            4 & 0 & 2 \\
            11 & 3 & 4 \\
            \end{array}\right]$,
    \item ${\bf a}^{(1)}=(\overline{2,1})$,
    \item ${\bf a}^{(2)}=(\overline{1,0})$,
    \item ${\bf b}^{(1)}=(0,2)$,
    \item ${\bf b}^{(2)}=(-1,1)$.
\end{itemize}   
\end{Example}

\begin{Example}\label{E3}
    Let $(x^{(1)},x^{(2)})$ and $(y^{(1)},y^{(2)})$ have the continued fraction expansions $[[-2,\overline{1,2}],[\overline{1,0}]]$ and $[[-3,\overline{1,0}],[\overline{0,1}]]$, respectively. Let us compute the beginning terms of the continued fraction expansion of $\left(x^{(1)}+y^{(2)},x^{(2)}y^{(1)}\right)$. Hence, our input is
    \begin{itemize}
        \item ${\bf a}^{(1)}=(-2,\overline{1,2})$,
        \item ${\bf a}^{(2)}=(\overline{1,0})$,
        \item ${\bf b}^{(1)}=(-3,\overline{1,0})$,
        \item ${\bf b}^{(2)}=(\overline{0,1})$,
        \item $C^{(1)}=\left[\begin{array}{ccc}
            0 & 0 & 1 \\
            0 & 0 & 0 \\
            0 & 1 & 0
        \end{array}\right]$,
        \item $C^{(2)}=\left[\begin{array}{ccc}
            0 & 0 & 0 \\
            1 & 0 & 0 \\
            0 & 0 & 0
        \end{array}\right]$,
        \item $C^{(3)}=\left[\begin{array}{ccc}
            0 & 0 & 0 \\
            0 & 0 & 0 \\
            0 & 0 & 1
        \end{array}\right]$,
    \end{itemize}
    whereas our output is $({\bf d}^{(1)},{\bf d}^{(2)})$.

    According to the algorithm, we start with input transformation, input swapping, and one more input transformation.
    \begin{align*}
        C^{(1)}\leftarrow &\left[\begin{array}{ccc}
            -2 & 1 & 1 \\
            1 & 0 & 0 \\
            0 & 1 & 0
        \end{array}\right]\cdot\left[\begin{array}{ccc}
            0 & 0 & 1 \\
            0 & 0 & 0 \\
            0 & 1 & 0
        \end{array}\right]=\left[\begin{array}{ccc}
            0 & 1 & -2 \\
            0 & 0 & 1 \\
            0 & 0 & 0 \\
            \end{array}\right]\\
        C^{(2)}\leftarrow &\left[\begin{array}{ccc}
            -2 & 1 & 1 \\
            1 & 0 & 0 \\
            0 & 1 & 0
        \end{array}\right]\cdot\left[\begin{array}{ccc}
            0 & 0 & 0 \\
            1 & 0 & 0 \\
            0 & 0 & 0
        \end{array}\right]=\left[\begin{array}{ccc}
            1 & 0 & 0 \\
            0 & 0 & 0 \\
            1 & 0 & 0 \\
            \end{array}\right]\\
        C^{(3)}\leftarrow &\left[\begin{array}{ccc}
            -2 & 1 & 1 \\
            1 & 0 & 0 \\
            0 & 1 & 0
        \end{array}\right]\cdot\left[\begin{array}{ccc}
            0 & 0 & 0 \\
            0 & 0 & 0 \\
            0 & 0 & 1
        \end{array}\right]=\left[\begin{array}{ccc}
            0 & 0 & 1 \\
            0 & 0 & 0 \\
            0 & 0 & 0 \\
            \end{array}\right]\\
        {\bf a}^{(1)}\leftarrow &(\overline{1,2}),\,
        {\bf a}^{(2)}\leftarrow (\overline{0,1})
    \end{align*}
    \begin{align*}
        C^{(1)}\leftarrow &\left[\begin{array}{ccc}
            0 & 0 & 0 \\
            1 & 0 & 0 \\
            -2 & 1 & 0 \\
            \end{array}\right]\\
        C^{(2)}\leftarrow &\left[\begin{array}{ccc}
            1 & 0 & 1 \\
            0 & 0 & 0 \\
            0 & 0 & 0 \\
            \end{array}\right]\\
        C^{(3)}\leftarrow &\left[\begin{array}{ccc}
            0 & 0 & 0 \\
            0 & 0 & 0 \\
            1 & 0 & 0 \\
            \end{array}\right]\\
        {\bf a}^{(1)}\leftarrow &(-3,\overline{1,0}),\,
        {\bf a}^{(2)}\leftarrow (\overline{0,1}),\,
        {\bf b}^{(1)}\leftarrow (\overline{1,2}),\,
        {\bf b}^{(2)}\leftarrow (\overline{0,1})
    \end{align*}
    \begin{align*}
        C^{(1)}\leftarrow &\left[\begin{array}{ccc}
            -3 & 0 & 1 \\
            1 & 0 & 0 \\
            0 & 1 & 0
        \end{array}\right]\cdot\left[\begin{array}{ccc}
            0 & 0 & 0 \\
            1 & 0 & 0 \\
            -2 & 1 & 0 \\
            \end{array}\right]=\left[\begin{array}{ccc}
            -2 & 1 & 0 \\
            0 & 0 & 0\\
            1 & 0 & 0\\
            \end{array}\right]\\
        C^{(2)}\leftarrow &\left[\begin{array}{ccc}
            -3 & 0 & 1 \\
            1 & 0 & 0 \\
            0 & 1 & 0
        \end{array}\right]\cdot\left[\begin{array}{ccc}
            1 & 0 & 1 \\
            0 & 0 & 0 \\
            0 & 0 & 0 \\
            \end{array}\right]=\left[\begin{array}{ccc}
            -3 & 0 & -3 \\
            1 & 0 & 1 \\
            0 & 0 & 0 \\
            \end{array}\right]\\
        C^{(3)}\leftarrow &\left[\begin{array}{ccc}
            -3 & 0 & 1 \\
            1 & 0 & 0 \\
            0 & 1 & 0
        \end{array}\right]\cdot\left[\begin{array}{ccc}
            0 & 0 & 0 \\
            0 & 0 & 0 \\
            1 & 0 & 0 \\
        \end{array}\right]=\left[\begin{array}{ccc}
            1 & 0 & 0 \\
            0 & 0 & 0 \\
            0 & 0 & 0 \\
            \end{array}\right]
            \\
        {\bf a}^{(1)}\leftarrow &(\overline{1,0}),\,
        {\bf a}^{(2)}\leftarrow (\overline{1,0})
    \end{align*}
    Next we perform input swaps and input transformations until we get the equalities $\left\lfloor\frac{c_{j_1,l_1}^{(1)}}{c_{j_1,l_1}^{(3)}}\right\rfloor =\left\lfloor\frac{c_{j_2,l_2}^{(1)}}{c_{j_2,l_2}^{(3)}}\right\rfloor$ and $\left\lfloor\frac{c_{j_1,l_1}^{(2)}}{c_{j_1,l_1}^{(3)}}\right\rfloor =\left\lfloor\frac{c_{j_2,l_2}^{(2)}}{c_{j_2,l_2}^{(3)}}\right\rfloor$, $j_1,j_2,l_1,l_2\in\{1,2,3\}$.
    \begin{align*}
        C^{(1)}\leftarrow &\left[\begin{array}{ccc}
            -2 & 0 & 1 \\
            1 & 0 & 0\\
            0 & 0 & 0\\
            \end{array}\right]\\
        C^{(2)}\leftarrow &\left[\begin{array}{ccc}
            -3 & 1 & 0 \\
            0 & 0 & 0 \\
            -3 & 1 & 0 \\
            \end{array}\right]\\
        C^{(3)}\leftarrow &\left[\begin{array}{ccc}
            1 & 0 & 0 \\
            0 & 0 & 0 \\
            0 & 0 & 0 \\
            \end{array}\right]\\
        {\bf a}^{(1)}\leftarrow &(\overline{1,2}),\,
        {\bf a}^{(2)}\leftarrow (\overline{0,1}),\,
        {\bf b}^{(1)}\leftarrow (\overline{1,0}),\,
        {\bf b}^{(2)}\leftarrow (\overline{1,0})
    \end{align*}
    \begin{align*}
        C^{(1)}\leftarrow &\left[\begin{array}{ccc}
            1 & 0 & 1 \\
            1 & 0 & 0 \\
            0 & 1 & 0
        \end{array}\right]\cdot\left[\begin{array}{ccc}
            -2 & 0 & 1 \\
            1 & 0 & 0\\
            0 & 0 & 0\\
            \end{array}\right]=\left[\begin{array}{ccc}
            -2 & 0 & 1 \\
            -2 & 0 & 1 \\
            1 & 0 & 0\\
            \end{array}\right]\\
        C^{(2)}\leftarrow &\left[\begin{array}{ccc}
            1 & 0 & 1 \\
            1 & 0 & 0 \\
            0 & 1 & 0
        \end{array}\right]\cdot\left[\begin{array}{ccc}
            -3 & 1 & 0 \\
            0 & 0 & 0 \\
            -3 & 1 & 0 \\
            \end{array}\right]=\left[\begin{array}{ccc}
            -6 & 2 & 0 \\
            -3 & 1 & 0 \\
            0 & 0 & 0 \\
            \end{array}\right]\\
        C^{(3)}\leftarrow &\left[\begin{array}{ccc}
            1 & 0 & 1 \\
            1 & 0 & 0 \\
            0 & 1 & 0
        \end{array}\right]\cdot\left[\begin{array}{ccc}
            1 & 0 & 0 \\
            0 & 0 & 0 \\
            0 & 0 & 0 \\
            \end{array}\right]=\left[\begin{array}{ccc}
            1 & 0 & 0 \\
            1 & 0 & 0 \\
            0 & 0 & 0 \\
            \end{array}\right]\\
        {\bf a}^{(1)}\leftarrow &(\overline{2,1}),\,
        {\bf a}^{(2)}\leftarrow (\overline{1,0})
    \end{align*}
    \begin{align*}
        C^{(1)}\leftarrow &\left[\begin{array}{ccc}
            -2 & -2 & 1 \\
            0 & 0 & 0 \\
            1 & 1 & 0\\
            \end{array}\right]\\
        C^{(2)}\leftarrow &\left[\begin{array}{ccc}
            -6 & -3 & 0 \\
            2 & 1 & 0 \\
            0 & 0 & 0 \\
            \end{array}\right]\\
        C^{(3)}\leftarrow &\left[\begin{array}{ccc}
            1 & 1 & 0 \\
            0 & 0 & 0 \\
            0 & 0 & 0 \\
            \end{array}\right]\\
        {\bf a}^{(1)}\leftarrow &(\overline{1,0}),\,
        {\bf a}^{(2)}\leftarrow (\overline{1,0}),\,
        {\bf b}^{(1)}\leftarrow (\overline{2,1}),\,
        {\bf b}^{(2)}\leftarrow (\overline{1,0})
    \end{align*}
    \begin{align*}
        C^{(1)}\leftarrow &\left[\begin{array}{ccc}
            1 & 1 & 1 \\
            1 & 0 & 0 \\
            0 & 1 & 0
        \end{array}\right]\cdot\left[\begin{array}{ccc}
            -2 & -2 & 1 \\
            0 & 0 & 0 \\
            1 & 1 & 0\\
            \end{array}\right]=\left[\begin{array}{ccc}
            -1 & -1 & 1 \\
            -2 & -2 & 1 \\
            0 & 0 & 0 \\
            \end{array}\right]\\
        C^{(2)}\leftarrow &\left[\begin{array}{ccc}
            1 & 1 & 1 \\
            1 & 0 & 0 \\
            0 & 1 & 0
        \end{array}\right]\cdot\left[\begin{array}{ccc}
            -6 & -3 & 0 \\
            2 & 1 & 0 \\
            0 & 0 & 0 \\
            \end{array}\right]=\left[\begin{array}{ccc}
            -4 & -2 & 0 \\
            -6 & -3 & 0 \\
            2 & 1 & 0 \\
            \end{array}\right]\\
        C^{(3)}\leftarrow &\left[\begin{array}{ccc}
            1 & 1 & 1 \\
            1 & 0 & 0 \\
            0 & 1 & 0
        \end{array}\right]\cdot\left[\begin{array}{ccc}
            1 & 1 & 0 \\
            0 & 0 & 0 \\
            0 & 0 & 0 \\
            \end{array}\right]=\left[\begin{array}{ccc}
            1 & 1 & 0 \\
            1 & 1 & 0 \\
            0 & 0 & 0 \\
            \end{array}\right]\\
        {\bf a}^{(1)}\leftarrow &(\overline{0,1}),\,
        {\bf a}^{(2)}\leftarrow (\overline{0,1})
    \end{align*}
    \begin{align*}
        C^{(1)}\leftarrow &\left[\begin{array}{ccc}
            -1 & -2 & 0 \\
            -1 & -2 & 0 \\
            1 & 1 & 0 \\
            \end{array}\right]\\
        C^{(2)}\leftarrow &\left[\begin{array}{ccc}
            -4 & -6 & 2 \\
            -2 & -3 & 1 \\
            0 & 0 & 0 \\
            \end{array}\right]\\
        C^{(3)}\leftarrow &\left[\begin{array}{ccc}
            1 & 1 & 0 \\
            1 & 1 & 0 \\
            0 & 0 & 0 \\
            \end{array}\right]\\
        {\bf a}^{(1)}\leftarrow &(\overline{2,1}),\,
        {\bf a}^{(2)}\leftarrow (\overline{1,0}),\,
        {\bf b}^{(1)}\leftarrow (\overline{0,1}),\,
        {\bf b}^{(2)}\leftarrow (\overline{0,1})
    \end{align*}
    \begin{align*}
        C^{(1)}\leftarrow &\left[\begin{array}{ccc}
            2 & 1 & 1 \\
            1 & 0 & 0 \\
            0 & 1 & 0
        \end{array}\right]\cdot\left[\begin{array}{ccc}
            -1 & -2 & 0 \\
            -1 & -2 & 0 \\
            1 & 1 & 0 \\
            \end{array}\right]=\left[\begin{array}{ccc}
            -2 & -5 & 0 \\
            -1 & -2 & 0 \\
            -1 & -2 & 0 \\
            \end{array}\right]\\
        C^{(2)}\leftarrow &\left[\begin{array}{ccc}
            2 & 1 & 1 \\
            1 & 0 & 0 \\
            0 & 1 & 0
        \end{array}\right]\cdot\left[\begin{array}{ccc}
            -4 & -6 & 2 \\
            -2 & -3 & 1 \\
            0 & 0 & 0 \\
            \end{array}\right]=\left[\begin{array}{ccc}
            -10 & -15 & 5 \\
            -4 & -6 & 2 \\
            -2 & -3 & 1 \\
            \end{array}\right]\\
        C^{(3)}\leftarrow &\left[\begin{array}{ccc}
            2 & 1 & 1 \\
            1 & 0 & 0 \\
            0 & 1 & 0
        \end{array}\right]\cdot\left[\begin{array}{ccc}
            1 & 1 & 0 \\
            1 & 1 & 0 \\
            0 & 0 & 0 \\
            \end{array}\right]=\left[\begin{array}{ccc}
            3 & 3 & 0 \\
            1 & 1 & 0 \\
            1 & 1 & 0 \\
            \end{array}\right]\\
        {\bf a}^{(1)}\leftarrow &(\overline{1,2}),\,
        {\bf a}^{(2)}\leftarrow (\overline{0,1})
    \end{align*}
    \begin{align*}
        C^{(1)}\leftarrow &\left[\begin{array}{ccc}
            -2 & -1 & -1 \\
            -5 & -2 & -2 \\
            0 & 0 & 0 \\
            \end{array}\right]\\
        C^{(2)}\leftarrow &\left[\begin{array}{ccc}
            -10 & -4 & -2 \\
            -15 & -6 & -3 \\
            5 & 2 & 1 \\
            \end{array}\right]\\
        C^{(3)}\leftarrow &\left[\begin{array}{ccc}
            3 & 1 & 1 \\
            3 & 1 & 1 \\
            0 & 0 & 0 \\
            \end{array}\right]\\
        {\bf a}^{(1)}\leftarrow &(\overline{0,1}),\,
        {\bf a}^{(2)}\leftarrow (\overline{0,1}),\,
        {\bf b}^{(1)}\leftarrow (\overline{1,2}),\,
        {\bf b}^{(2)}\leftarrow (\overline{0,1})
    \end{align*}
    \begin{align*}
        C^{(1)}\leftarrow &\left[\begin{array}{ccc}
            0 & 0 & 1 \\
            1 & 0 & 0 \\
            0 & 1 & 0
        \end{array}\right]\cdot\left[\begin{array}{ccc}
            -2 & -1 & -1 \\
            -5 & -2 & -2 \\
            0 & 0 & 0 \\
            \end{array}\right]=\left[\begin{array}{ccc}
            0 & 0 & 0 \\
            -2 & -1 & -1 \\
            -5 & -2 & -2 \\
            \end{array}\right]\\
        C^{(2)}\leftarrow &\left[\begin{array}{ccc}
            0 & 0 & 1 \\
            1 & 0 & 0 \\
            0 & 1 & 0
        \end{array}\right]\cdot\left[\begin{array}{ccc}
            -10 & -4 & -2 \\
            -15 & -6 & -3 \\
            5 & 2 & 1 \\
            \end{array}\right]=\left[\begin{array}{ccc}
            5 & 2 & 1 \\
            -10 & -4 & -2 \\
            -15 & -6 & -3 \\
            \end{array}\right]\\
        C^{(3)}\leftarrow &\left[\begin{array}{ccc}
            0 & 0 & 1 \\
            1 & 0 & 0 \\
            0 & 1 & 0
        \end{array}\right]\cdot\left[\begin{array}{ccc}
            3 & 1 & 1 \\
            3 & 1 & 1 \\
            0 & 0 & 0 \\
            \end{array}\right]=\left[\begin{array}{ccc}
            0 & 0 & 0 \\
            3 & 1 & 1 \\
            3 & 1 & 1 \\
            \end{array}\right]\\
        {\bf a}^{(1)}\leftarrow &(\overline{1,0}),\,
        {\bf a}^{(2)}\leftarrow (\overline{1,0})
    \end{align*}

    Proceeding in this way, after $20$ inputs transformations, we arrive to the following transformation matrices:

     \begin{align*}
        C^{(1)}\leftarrow &\left[\begin{array}{ccc}
            -2707 & -2149 & -705 \\
            -718 & -570 & -187 \\
            -1436 & -1140 & -374 \\
            \end{array}\right]\\
        C^{(2)}\leftarrow &\left[\begin{array}{ccc}
            -7080 & -5260 & -1840 \\
            -1770 & -1405 & -460 \\
            -3894 & -3091 & -1012 \\
            \end{array}\right]\\
        C^{(3)}\leftarrow &\left[\begin{array}{ccc}
            2453 & 1947 & 638 \\
            669  & 531 &  174\\
            1338  & 1062 & 348 \\
            \end{array}\right]
    \end{align*}
We are now able to extract the first output of the MCF of $\left(x^{(1)}+y^{(2)},x^{(2)}y^{(1)}\right)$. In fact, it is possible to see that for all $j_1,j_2,l_1,l_2\in\{1,2,3\}$ we have
\begin{align*}
\left\lfloor\frac{c_{j_1,l_1}^{(1)}}{c_{j_1,l_1}^{(3)}}\right\rfloor &=\left\lfloor\frac{c_{j_2,l_2}^{(1)}}{c_{j_2,l_2}^{(3)}}\right\rfloor=-2,\\
\left\lfloor\frac{c_{j_1,l_1}^{(2)}}{c_{j_1,l_1}^{(3)}}\right\rfloor &=\left\lfloor\frac{c_{j_2,l_2}^{(2)}}{c_{j_2,l_2}^{(3)}}\right\rfloor=-3.
\end{align*}
After the first output, the sequence of partial quotients is $[[-2,\ldots],[-3,\ldots]]$, and we can perform the output transformation:

    \begin{align*}
        C^{(1)}&\leftarrow \left[\begin{array}{ccc}
            2453 & 669 & 1338 \\
            1947 & 531 & 1062 \\
            638 & 174 & 348 \\
            \end{array}\right]\\
        C^{(2)}&\leftarrow\left[\begin{array}{ccc}
            -2707 & -718 & -1436 \\
            -2149 & -570 & -1140 \\
            -705 & -187 & -374 \\
            \end{array}\right]+2\left[\begin{array}{ccc}
            2453 & 669 & 1338 \\
            1947 & 531 & 1062 \\
            638 & 174 & 348 \\
            \end{array}\right]=\left[\begin{array}{ccc}
            2199 & 620 & 1240 \\
            1745 & 492 & 984 \\
            571 & 161 & 322 \\
            \end{array}\right]\\
        C^{(3)}&\leftarrow\left[\begin{array}{ccc}
            -7080 & -1770 & -3894 \\
            -5620 & -1405 & -3091 \\
            -1840 & -460 & -1012 \\
            \end{array}\right]+3\left[\begin{array}{ccc}
            2453 & 669 & 1338 \\
            1947 & 531 & 1062 \\
            638 & 174 & 348 \\
            \end{array}\right]=\left[\begin{array}{ccc}
            279 & 237 & 120 \\
            221 & 188 & 95 \\
        74 & 62 & 32 \\
            \end{array}\right]
    \end{align*}

By keeping performing the algorithm using the latter matrices,  we get that the continued fraction expansion for $\left(x^{(1)}+y^{(2)},x^{(2)}y^{(1)}\right)$ is $$[[-2, 7, 36, 1, 3, 4, 2,...], [-3, 7, 26, 0, 0, 0, 1,...]].$$
\end{Example}


\begin{thebibliography}{100}
\bibitem{BGSch} M. Beeler, R. W. Gosper, R. Schr\"oppel, \textit{HAKMEM, A. I. Lab Memo \# 239}, M.I.T. 1972.

\bibitem{Ber1} L. Bernstein, \textit{Periodical continued fractions for irrationals n by Jacobi’s algorithm}.
\newblock J. Reine Angew. Math. 213 (1964), 31–38.

\bibitem{Ber2} L. Bernstein, \textit{Periodicity of Jacobi’s algorithm for a special case of cubic irrationals}.
\newblock J. Reine Angew. Math. 213 (1964), 137–147.

\bibitem{Ber3} L. Bernstein, \textit{New infinite classes of periodic Jacobi–Perron algorithms}. 
\newblock Pacific J. Math. 16(3) (1965), 439–469.

\bibitem{Berthe} V. Berthé, K. Dajani, C. Kalle, E. Krawczyk, H. Kuru, A. Thevis, \textit{Rational Approximations, Multidimensional Continued Fractions, and Lattice Reduction}. 
\newblock In: Abdellatif, R., Karemaker, V., Smajlovic, L. (eds) Women in Numbers Europe IV, Association for Women in Mathematics Series, vol 32. Springer, Cham (2024).

\bibitem{Col}{M. J. Collins},
\textit{Arithmetic on Continued Fractions.},
\newblock preprint, (2024), arXiv: \href{https://arxiv.org/abs/2412.19929}{2412.19929}.

\bibitem{Cusick} T. W. Cusick, \textit{Sums and products of continued fractions}.
\newblock Proceedings of the American Mathematical Society 27(1) (1971), 35--38.

\bibitem{Dubois} E. Dubois, A. Farhane, R. Paysant-Le Roux, \textit{Etude des interruptions dans l'algorithme de Jacobi-Perron}.
\newblock Bulletin of the Australian Mathematical Society 69 (2004), 241-254.


\bibitem{Gosper} R.W. Gosper, \textit{Continued Fraction Arithmetic}, preprint 1976.

\bibitem{Hall} M. Hall, \textit{On the Sum and Product of Continued Fractions}. 
\newblock Annals of Mathematics 48(4) (1947), 966-993.

\bibitem{HER}{C. Hermite},
\textit{Extraits de lettres de M. Ch. Hermite à M. Jacobi sur différents objects de la théorie des nombres},
\newblock J. Reine Angew. Math. \textbf{40} (1850), 261-278.

\bibitem{Hur}{A. Hurwitz},
\textit{Über die Kettenbruch-Entwicklung der Zahl e},
\newblock Schriften der physikalisch-ökonomischen Gesellschaft zu Königsberg in Pr., 32. Jahrgang, (1891). Mathematische Werke: Zweiter Band Zahlentheorie Algebra und Geometrie (1963), 129-133.




\bibitem{Jacobi} C. G. J. Jacobi, \textit{Correspondance mathematique avec Legendre.} 
\newblock In Gesammelte Werke 6 (Berlin Academy, Berlin, 1891), 385-426.

\bibitem{Lev} C. Levesque, G. Rhin, \textit{Two families of periodic Jacobi algorithms with period lengths going to infinity}.
\newblock J. Number Theory 37 (1991), 173-180.

\bibitem{LiarStam} P. Liardet, P. Stambul, \textit{Algebraic Computations with Continued Fractions}.
\newblock J. Number Theory 73 (1998), 92-121.

\bibitem{Perron} O. Perron, \textit{Grundlagen f\"ur eine theorie des Jacobischen kettenbruchalgorithmus.}
\newblock Math. Ann. 64 (1907), 1–76.

\bibitem{Rada} H. \v{R}ada, S. Starosta, V. Kala, \textit{Periodicity of general multidimensional continued fractions using repetend matrix form}.
\newblock Expositiones Mathematicae 42(3) (2024), Article ID 125571.

\bibitem{Raney} G. Raney, \textit{On Continued Fractions and Finite Automata}
\newblock Mathematische Annalen 206 (1973), 265-283.

\bibitem{Voutier} P. Voutier, \textit{Families of periodic Jacobi–Perron algorithms for all period lengths.}
\newblock J. Number Theory 168 (2016), 472-486.
\end{thebibliography}
\end{document}